 \documentclass[final,3p,times,12pt]{elsarticle}

\usepackage{amssymb}
\usepackage{graphicx,wrapfig,lipsum}
\usepackage{amsmath}
\usepackage{tabularx}
\usepackage{adjustbox}
\usepackage{amssymb}
\usepackage{caption}
\usepackage{subcaption}
\usepackage{fancyref}
\usepackage{float}
\usepackage{adjustbox}
\usepackage{multirow}
\usepackage{enumerate}
\usepackage{amsthm}
\usepackage[utf8]{inputenc}
\usepackage{yfonts}
\usepackage{mathtools}
\usepackage{algorithm}
\usepackage{algpseudocode}
\usepackage{afterpage}

\newtheorem{remark}{Remark}

\newcommand{\norm}[1]{\left\lVert#1\right\rVert}
\usepackage{xcolor}
\usepackage{array}
\usepackage{tabu}

\newtheorem*{acknowledgments}{Acknowledgments}
\newtheorem*{CoreIdeas}{Core Ideas}

\newcommand{\uw}{{\bf{u_w}}}
\newcommand{\uwn}{{\bf{u^n_w}}}
\newcommand{\uwo}{{\bf{u^{n,j}_w}}}
\newcommand{\uwj}{{\bf{u^{n,j+1}_w}}}
\newcommand{\ez}{{\bf{e_z}}}

\DeclareMathOperator{\sign}{sign}

\usepackage{lineno}

\begin{document}

\begin{frontmatter}



\title{Efficient Solvers for Nonstandard Models for Flow and Transport in
	Unsaturated Porous Media} 
\author[add1]{Davide Illiano}
\ead{Davide.Illiano@uin.no}
\author[add1]{Jakub Wiktor Both}
\ead{Jakub.Both@uib.no}
\author[add1,add2]{Iuliu Sorin Pop}
\ead{sorin.pop@uhasselt.be}
\author[add1]{Florin Adrian Radu}
\ead{Florin.Radu@uib.no}

\address[add1]{Department of Mathematics, University of Bergen, Allegaten 41, Bergen, Norway}
\address[add2]{Faculty of Science, University of Hasselt,Agoralaan Building D, BE 3590 Diepenbeek, Belgium}

\begin{abstract}
	We study several iterative methods for fully coupled flow and reactive transport in porous media. The resulting mathematical model is a coupled, nonlinear evolution system. The flow model component builds on the Richards equation, modified to incorporate nonstandard effects like dynamic capillarity and hysteresis, and a reactive transport equation for the solute. 
The two model components are strongly coupled. On one hand, the flow affects the concentration of the solute; on the other hand, the surface tension is a function of the solute, which impacts the capillary pressure and, consequently, the flow.
After applying an Euler implicit scheme, we consider a set of iterative linearization schemes to solve the resulting nonlinear equations, including both monolithic and two splitting strategies. The latter include a canonical nonlinear splitting and an alternate linearized splitting, which appears to be overall faster in terms of numbers of iterations, based on our numerical studies.
The (time discrete) system being nonlinear, we investigate different linearization methods. We consider the linearly convergent L-scheme, which converges unconditionally, and the Newton method, converging quadratically but subject to restrictions on the initial guess. Whenever hysteresis effects are included, the Newton method fails to converge. The L-scheme converges; nevertheless, it may require many iterations. This aspect is improved by using the Anderson acceleration.
A thorough comparison of the different solving strategies is presented in five numerical examples, implemented in MRST, a toolbox based on MATLAB.

\end{abstract}

\end{frontmatter}


\section{Introduction}

Mathematical models for complex physical phenomena are generally neglecting several processes, in order to guarantee that the result is sufficiently simple and to facilitate the numerical simulations. With a particular focus on porous media applications, in this sense we mention enhanced oil recovery, diffusion of substances in living tissues, and pollution of underground aquifers. With the increase of computational power, and the development of efficient simulation algorithms, mathematical models are improved continuously, and more and more of the neglected effects are included.

\par 
When studying unsaturated flow, the equilibrium capillary pressure plays a fundamental role. It is typically assumed to be a nonlinear, monotone function of the water content. Explicit representations have been obtained thanks to numerous experiments under equilibrium conditions (no flowing phases). Even though this formulation has been the most commonly used in the last decades, it has been observed \cite{Camps-Roach,DiCarlo,Oung,stauffer1978}, that changes in time of the water content, thus its time derivatives, do influence the profile of the capillary pressure. In terms of modeling, this is achieved by including the so-called dynamic effects \cite{Beliaev2001, Gray, Mikelic2009}. Numerous papers investigate the existence of a solution for systems including such effects, among them we cite \cite{Cao,Koch,Milisic}. Furthermore, the problem has been already studied numerically in, e.g.,  \cite{Abreu2017,Abreu2020,Cao2019}.

\par
\begin{figure}[h]
\centering
\includegraphics[scale=.55]{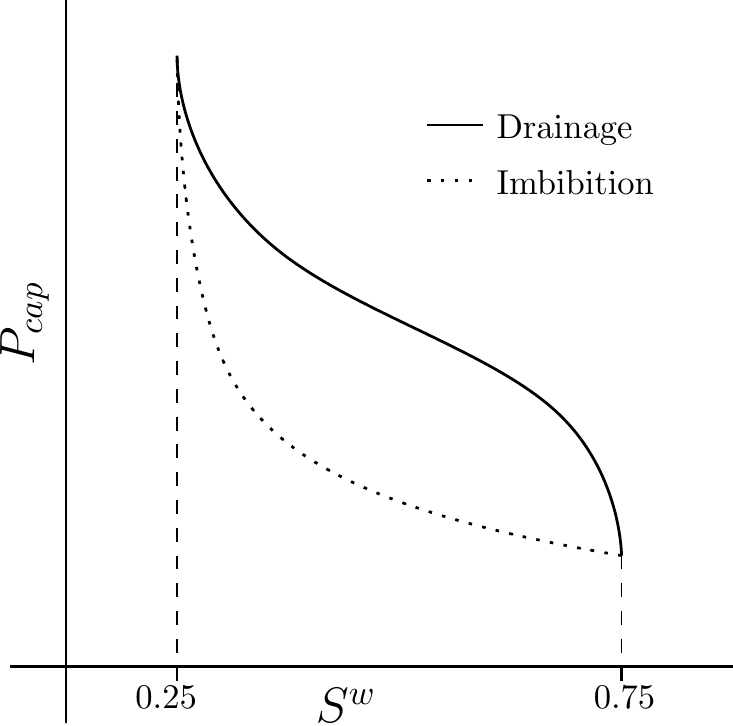}
\caption{Primary hysteresis loop as presented in \cite{McClure2018}.}
\label{HysteresisLoop}
\end{figure}

The hysteresis effect is another phenomenon often neglected. Again, experiments have revealed that the curve obtained when investigating the imbibition process, is different from the one observed during the drainage, \cite{Carroll,Hoa1977,McClure2018,Morrow}. This is sketched in Fig. \ref{HysteresisLoop}.

\par 
In this article, we study unsaturated flow in porous media, modeled by the Richards equation \cite{RichardsBerardi,BookHelmig}, however including both dynamic and hysteresis effects. Furthermore, we include a solute component, e.g., a surfactant, in the wetting phase, which can directly influence the fluid properties (\cite{agosti,prechtel}). The study of the transportation of an external components, e.g., surfactant, in variably saturated porous media has been already investigated both numerically \cite{Illiano2020,surf,Radu2013} and experimentally \cite{surfactant1999,Karagunduz}. 
\par
Here, we will mainly concentrate on numerical
studies, extending the solution techniques in  \cite{Illiano2020} to include dynamic capillarity and hysteresis
We consider the following model for coupled unsaturated flow and reactive transport
\begin{equation}\label{standard}
\begin{split}
\partial_t \theta(\Psi,c) - \nabla \cdot (K(\theta(\Psi,c))\nabla(\Psi+z))\ =\ \mathbb{S}_1, \\
\partial_t (\theta(\Psi,c) c) - \nabla \cdot (D\nabla c - \uw c) + R(c)\ =\ \mathbb{S}_2.\\
\end{split}
\end{equation}
Here, $\theta(\Psi,c)$ is the water content, expressed as a function of both the unknown pressure head $\Psi$ and the concentration of the external component $c$. $K$, a function of the water content $\theta$, is the conductivity, $z$ the vertical coordinate of $\vec{x}$, pointing against gravity, $D$ the dispersion/diffusion coefficient, $\uw:=-K(\theta(\Psi,c))\nabla(\Psi+z)$ the water flux, $R(c)$ the reaction term and $\mathbb{S}_1,\mathbb{S}_2$ the external sink/source terms involved.
\par

Next to a concentration dependence of $ \theta $ and $ \Psi $, here we include also play-type hysteresis and dynamic capillary effects as introduced in \cite{Beliaev2001}. More precisely,
\begin{equation}\label{nonstandard}
\Psi \in - p_{cap}(\theta,c) + \tau(\theta)\partial_t \theta  + \gamma(\theta) \sign (\partial_t\theta) ,
\end{equation}
where $p_{cap}$ is the equilibrium capillary pressure, expressed as a function of $\theta$ and $c$, $\tau(\theta)$ the dynamic effects, and $\gamma(\theta)$ the width of the primary hysteresis loop. Later on, for ease of presentation, we consider $\gamma$ as a positive constant, $\gamma\in\mathbb{R}_{\geq 0}$. Note that (\ref{nonstandard}) is a differential inclusion as the sign graph is multi-valued and defined as follow,
\begin{equation}
	\begin{split}
		sign(\xi) = \begin{cases} 1 &\mbox{for}\  \xi>0,\\
			[-1,1]  &\mbox{for}\ \xi=0,\\
			-1 &\mbox{for}\ \xi<0. \end{cases}
	\end{split}
\end{equation}
The multi-valued graph allows switching between the imbibition and drainage curves in the play-type hysteresis. For more details on the formulation we refer to \cite{Beliaev2001}. 

The primary unknowns of the system are the pressure $\Psi$, the concentration $c$ and the water content $\theta$. In standard models, also obtained as special case for $\gamma=\tau(\theta)=0$, $\theta$ is a function of pressure and concentration. Therefore, (\ref{nonstandard}) is replaced by an algebraic relationship, which simplifies the model and allows eliminating $ \theta $ as an unknown.
  In the extended/nonstandard formulation, $\theta$ is an unknown and (\ref{nonstandard}) is required as additional equation of the model. Initial and boundary conditions will complete the system. 

To avoid working with a graph, we consider the following regularization,
\begin{equation}\label{regularization}
\Phi(\xi) = \begin{cases} \sign(\xi) &\mbox{if}\qquad \vert \xi\vert \geq \delta, \\
\frac{\xi}{\delta} &\mbox{if}\qquad \vert \xi\vert < \delta,
\end{cases}
\end{equation}
where $\delta\in\mathbb{R}^{+}$ is a small parameter. Using this in (\ref{nonstandard}) gives the regularized system of equations
\begin{equation}\label{nonstandard_regularized}
\begin{split}
\partial_t \theta - \nabla \cdot (K(\theta)\nabla(\Psi+z))\ &=\ \mathbb{S}_1, \\
\Psi = - p_{cap}(\theta,c) + \tau(\theta)\partial_t \theta &+ \gamma \Phi(\partial_t\theta),\\
\partial_t (\theta c) - \nabla \cdot (D\nabla c - \uw c) + R(c)\ &=\ \mathbb{S}_2.\\
\end{split}
\end{equation}
From now on, the system \eqref{nonstandard_regularized} will be further investigated. We will discretize the equations and study different solving algorithms.
\begin{remark}
An inverse formulation is proposed in \cite{Beliaev2001}, obtained by solving \eqref{nonstandard}, as its regularized counterpart in \eqref{nonstandard_regularized}, in terms of $\partial_t \theta$. This gives
\begin{equation}\label{Ffunction}
 \partial_t \theta = F(\Psi,\theta,c),
\end{equation} 
for a suitable function $F$. The time derivative in the flow equation can then be substituted by $F$,
\begin{equation}
F(\Psi, \theta,c) - \nabla \cdot  \big(K(\theta)\nabla(\Psi+z)\big)\ =\ \mathbb{S}_1.
\end{equation}  
This formulation is used for the mathematical analysis of such models, \\\cite{Cao2015,Schweizer2012}.
It has been observed, e.g., in \cite{Lunowa2020}, that such formulation can reduce the number of iterations required to solve the system of equations, compared to the formulation in \eqref{nonstandard}.
However, for the particular test cases investigated here, no remarkable improvements are observed. Thus, for ease of presentation, we will report the results obtained only for the formulation given by (\ref{nonstandard_regularized}).
\end{remark}

\par
We point out that the concentration of the external component directly influences the capillary pressure. The presence of such a component results in a non-constant surface tension, which induces a rescaling of the pressures \cite{salt,smith1994,laboratory}. 
\par
To solve the system (\ref{nonstandard_regularized}) numerically, one first needs to discretize in time and space, and then develop solvers for the discretized equations.  In this paper, due to the expected low regularity of the solutions \cite{AltLuckhaus} and the desire of relatively large time steps, we choose to use the backward Euler method for the time discretization. Certain processes investigated in porous media flow can take place on time intervals longer than decades, thus the need for large time steps. Multiple spatial discretization techniques are available, e.g., the Galerkin Finite Element Method (\emph{FEM}) \cite{barrett1997,Nochetto,russell1983}, Discontinuous Galerkin Method (\emph{DGM}) \cite{Arnold2006,Karpinski,Li2007,Sun}, the Mixed Finite Element Method (\emph{MFEM}) \cite{arbogast1996,Cao2019,radu2010,Radu2013,vohralik,woodward}, the Finite Volume Method (\emph{FVM}) \\\cite{eymard1999} and the Multi-Point Flux Approximation (\emph{MPFA}) \\\citep{aavatsmark,Arraras2020,bause2010,Klausen}. We will here concentrate on \emph{FEM} and \emph{TPFA} (Two Points Flux Approximation), a particular case of \emph{MPFA}. In particular, we cite \cite{Berardi2020,Dolejsi,Zha} for papers on improved numerical schemes applied to the Richards equation.

\par 
Since the equations investigated here are characterized by several nonlinear quantities, $K(\theta)$, $p_{cap}(\theta,c)$, $\tau(\partial_t \theta)$, and $R(c)$, and the time discretization is not explicit, one needs to solve a nonlinear system at each time step, requiring a linearization procedure. Examples of possible linearization schemes are: the Newton method \cite{Paniconi}, the modified Picard method \cite{celia} and the L-scheme \cite{,List2016,Pop2004}. In this paper, we investigate the Newton method and the L-scheme. The former is a commonly used linearization scheme which is quadratically convergent. However, this convergence is only local and one needs to compute the Jacobian matrix, which can be expensive. The L-scheme is instead globally (linearly) convergent, under mild restrictions, and it does not require the computation of any derivative. The L-scheme is in general slower in terms of numbers of iterations than the Newton method. Moreover, the linear systems to be solved within each iteration are better conditioned when compared to the ones given by the Newton method \cite{Illiano2020,List2016}. Furthermore, the rate of convergence of the scheme strongly depends on user-defined parameters. Such aspects are investigated for numerous nonlinear problems, including Richards equation, and two-phase flow in porous media, in  \cite{Illiano2020,List2016,Mitra,Pop2004,Slodicka}. Finally there numerous papers proposing improved formulation of the L-scheme, among them we cite \cite{Albuja,Mitra}.
\par
 In this work, we test the L-scheme on more complex problems involving hysteresis and dynamic effects, and coupled reactive transport and flow. 
Furthermore, we investigate a post-processing technique, the Anderson Acceleration (AA) \cite{Anderson1965}, which can drastically improve linearly convergent schemes. The acceleration tool requires user-defined parameters. As will be seen below, choosing the suitable parameters for the AA, significantly relaxes the choice of the parameters for the L-scheme linearization. 
\par
We observe that the system (\ref{nonstandard_regularized}) is fully coupled. This is due to the dependence of the capillary pressure on both $\theta$ and $c$. Therefore, we will investigate multiple solution algorithms, combining different linearization schemes and decoupling techniques. Decoupling/splitting the equations may present multiple advantages such as: an easier implementation, a better conditioned problem to solve, similar convergence properties but faster computations. We divide the schemes into three main categories: monolithic (Mono), nonlinear splitting (NonLinS) and alternate splitting (AltS). Subsequently, we denote, e.g., by Newton-Mono, the monolithic scheme obtained by applying the Newton method as linearization. Such schemes have already been investigated for the standard model in \cite{Illiano2020}.

\par
The paper is organized as follows. In Section \ref{numericalrichards}, we present the linearization and discretization techniques including monolithic or decoupled solution approaches.  Section \ref{numerical} presents five different numerical examples, which allow to compare the efficiency and robustness of the solving algorithms.
Section \ref{conclusion} concludes this work with the final remarks.

\section{Problem formulation, discretization and iterative schemes}\label{numericalrichards}

In the following, we use the standard notations of functional analysis. The domain $\Omega\subset \mathbb{R}^d$, $d = 1,2$ or $3$, is bounded and has a Lipschitz continuous boundary $\partial \Omega$. The final time is $T>0$, and the time domain is $(0,T]$. $L^2(\Omega)$ denotes the space of real valued, square integrable functions defined on $\Omega$ and  $H^1(\Omega)$ its subspace containing the functions also having weak first derivatives in $L^2(\Omega)$. $H_0^1(\Omega)$ is the space of functions belonging to $H^1(\Omega)$, having zero trace on the boundary $\partial \Omega$. Furthermore, we denote by $< \cdot, \cdot > $ the standard $L^2(\Omega)$ scalar product and by $\norm{\cdot}$ the associated norm.

To numerically solve the system of equations (\ref{nonstandard_regularized}), one needs to discretize both in time and space. We combine the backward Euler method with linear Galerkin finite elements. Let $N \in \mathbb{N}$ be a strictly positive natural number. We define the time step size $\Delta t \ = \ T/N$ and $t_n\ =\ n\Delta t\ (n\in{1, 2, \dots, N})$. Furthermore, let $T_h$ be a regular decomposition of $\Omega$, $\overline{\Omega} = \underset{T \in T_h}{\cup} T$, with $h$ denoting the mesh diameter.  The finite element spaces $V_h \subset H_0^1(\Omega)$ and $W_h \subset L^2(\Omega)$ are defined by
\begin{equation}
V_h := \left\{ v_h\in H_0^1(\Omega)\ s.t.\ v_{h|T} \in \mathbb{P}_1(T), \ T\in T_h\right\}, \qquad W_h := \left\{ w_h\in L^2(\Omega)\ s.t.\ w_{h|T} \in \mathbb{P}_1(T), \ T\in T_h\right\},
\end{equation}
where $\mathbb{P}_1(T)$ denotes the space of the linear polynomials on $T$.
The fully discrete Galerkin formulation of the system (\ref{nonstandard_regularized}) can now be written as:

\textbf{Problem Pn:} Let $n \ge 1$ be fixed. Assuming that $\Psi^{n-1}_h, c^{n-1}_h \in V_h$ and $\theta^{n-1}_h \in W_h$ are given, find $\Psi^{n}_h, c^{n}_h \in V_h$ and $\theta^{n}_h \in W_h$ such that
\begin{equation}\label{richards1}
	\resizebox{0.92\hsize}{!}{$
\begin{split}
<\theta^{n}_h - \theta^{n-1}_h, v_{1,h} > + \Delta t <K(\theta^{n}_h)(\nabla\Psi^{n}_h+ {\ez}), \nabla v_{1,h} > \ =\ \Delta t <\mathbb{S}_1, v_{1,h} >\\
\Delta t<\Psi^{n}_h , w_{1,h}> +\Delta t <p_{cap}(\theta^{n}_h,c^{n}_h), w_{1,h}> - <\tau(\theta^{n}_h) (\theta^{n}_h - \theta^{n-1}_h), w_{1,h}> \ = \ \Delta t \gamma <\Phi\left(\frac{\theta^{n}_h - \theta^{n-1}_h}{\Delta t}\right), w_{1,h}>\\
< \theta^{n}_h (c^{n}_h  - c^{n-1}_h) +  c^{n}_h (\theta^{n}_h - \theta^{n-1}_h), v_{2,h} >+ \Delta t  <D \nabla \Psi^{n}_h + \uwn c^{n}_h, \nabla v_{2,h}>  + \Delta t <R(c^{n}_h),v_{2,h} > \ = \ \Delta t <\mathbb{S}_2, v_{2,h}>
\end{split}
$}
\end{equation}
holds for all $v_{1,h},v_{2,h} \in V_h$ and for all $w_{1,h}\in W_h$. We denote by ${\ez}$ the unit vector in the direction opposite to gravity.

Observe that choosing the space $ H^1_0(\Omega) $ implies that homogeneous boundary conditions have been adopted for the pressure and the concentration. However, this choice is made for the ease of presentation, the extension to other boundary conditions being possible without major complications. We also mention that, for $ n=1 $, we use the approximation in $ V_h $ of the initial water content and concentration, respectively $ \theta^0_h $ and $ c^0_h $.

In the following, we investigate different iterative schemes for solving Problem Pn. These schemes are based on the ones discussed in \cite{Illiano2020}, extending them, not only to the case of dynamic capillary pressure $(\tau(\theta)\neq 0)$ \cite{IllianoEnumath}, but also to the case of hysteresis. Among the numerous papers investigating numerically the effects of hysteresis and dynamic capillarity pressure, we cite \cite{Peszynska,Zhang}. As mentioned, we compare monolithic (Mono) and splitting (NonLinS and AltS) solvers, combined with two different linearization schemes, the Newton method and the L-scheme. Furthermore, the Anderson acceleration \cite{Anderson1965} will be taken into account to speed up the linearly convergent L-scheme.

\subsection{Solving algorithms}\label{section:algorithms}

In what follows, when solving \eqref{richards1} iteratively, the index $n$ will always refer to the time step level, whereas $j$ will denote the iteration index. As a rule, the iterations will start with the solution at the last time step, $t_{n-1}$, for example $\Psi^{n,1} = \Psi^{n-1}$. As mentioned, this choice is not required for L-type schemes but it is a natural one.

In a compact form Problem Pn can be seen as the system
\begin{equation}\label{system1}
\begin{cases} F_1 (\Psi_h^{n},  \theta_h^{n}) &= 0,\\
 F_2 (\Psi_h^{n}, \theta_h^{n},c_h^{n}) &= 0,\\
 F_3 (\Psi_h^{n}, \theta_h^{n},c_h^{n}) &= 0,
\end{cases}
\end{equation} 
with $F_1,F_2$ resulting from the flow equations and $F_3$ from the transport. In the following we will indicate with $F^{lin}$, the linearized formulation of $F$ obtained by either the Newton method or the L-scheme. Finally, we can proceed to present monolithic and splitting solvers.

\begin{figure}[h]
\centering
\includegraphics[scale=.6]{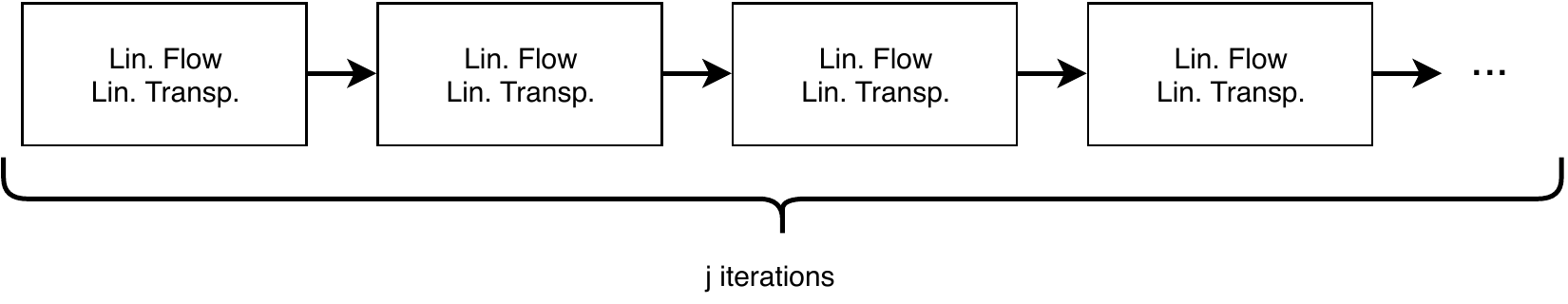}
\caption{The monolithic approach.}
\label{scheme0}
\end{figure}

In the monolithic approach one solves the three equations of the system \eqref{system1} at once. Formally, one iteration is:
\\Find  $\Psi_h^{n, j+1}$,  $\theta_h^{n, j+1}$ and $c_h^{n, j+1}$ such that
\begin{equation}
\begin{cases} F^{lin}_1 (\Psi_h^{n, j+1},  \theta_h^{n, j+1}) &= 0,\\
 F^{lin}_2 (\Psi_h^{n, j+1}, \theta_h^{n, j+1},c_h^{n, j+1}) &= 0,\\
 F^{lin}_3 (\Psi_h^{n, j+1}, \theta_h^{n, j+1},c_h^{n, j+1}) &= 0,
\end{cases}
\end{equation} 
where $F_i^{Lin}$ is the linearization of $F_i$, $i \in \{1,2,3\}$. Depending on which linearization technique is used, we refer to the Newton-monolithic scheme (Newton-Mono) or monolithic-$L$-scheme (LS-Mono). These two schemes will be presented in details below.
Fig. \ref{scheme0} displays the sketched version of the monolithic solver. 

\begin{figure}[h]
\centering
\includegraphics[scale=.55]{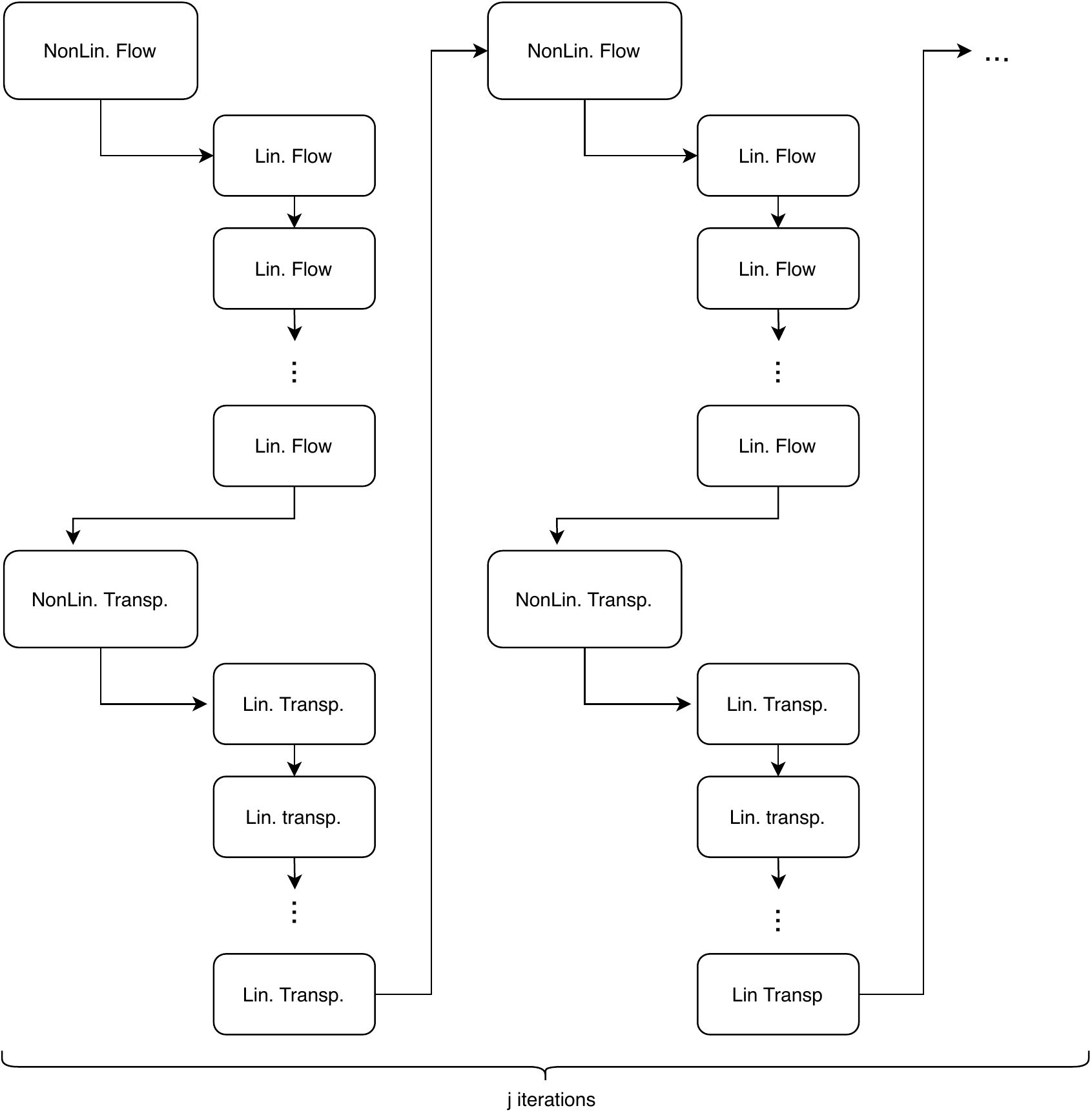}
\caption{The nonlinear splitting approach.}
\label{scheme2}
\end{figure}

\begin{figure}[h]
\centering
\includegraphics[scale=.6]{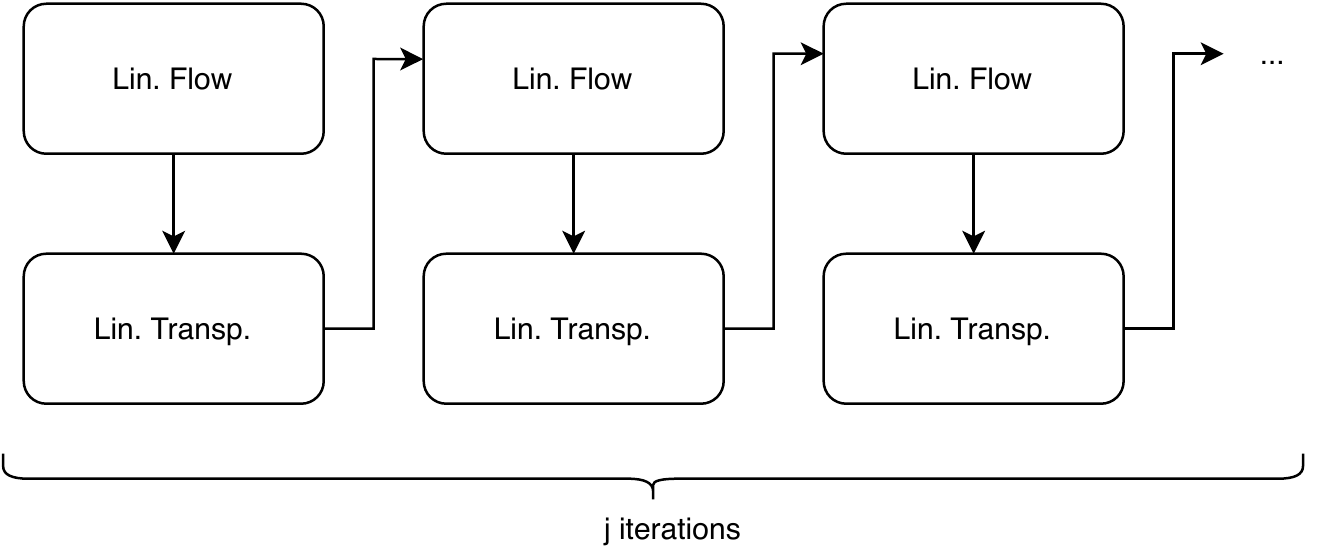}
\caption{The alternate splitting approach.}
\label{scheme3}
\end{figure}

In the iterative splitting approach, the flow and the transport equations are solved subsequently, iterating between them. We will distinguish between two primary splitting schemes: the nonlinear splitting (NonLinS) and the alternate linearized splitting (AltLinS), illustrated in Figure \ref{scheme2} and Figure \ref{scheme3}, respectively. Such schemes have already been studied, in the case of a standard flow model in \cite{Illiano2020}. 
\\In the nonlinear splitting, one iteration step is:
\\Find  first $\Psi_h^{n, j+1}$, $\theta_h^{n, j+1}$ such that
\begin{equation}\label{solving_system1}
\begin{split}
\begin{cases} F_1 (\Psi_h^{n, j+1}, \theta_h^{n,j+1}) &= 0,\\
 F_2 (\Psi_h^{n, j+1}, \theta_h^{n, j+1},c_h^{n, j}) &= 0,
 \end{cases}\\
 \end{split}
\end{equation}
and then find $c_h^{n,j+1}$ such that
\begin{equation}\label{eq13}
\begin{split}
F_3 (\Psi_h^{n, j+1},\theta_h^{n, j+1}, c_h^{n, j+1}) = 0.
\end{split}
\end{equation}
The two flow equations are solved at once. Each of the nonlinear systems \eqref{solving_system1} and \eqref{eq13} is solved until some convergence criterion is met. Once the pressure and water content are obtained, $\Psi_h^{n, j+1}$ and $ \theta_h^{n,j+1}$, are then used in the transport equation \eqref{eq13} to compute $c_h^{n, j+1}$.
The resulting $F_1$, $F_2$ and $F_3$, being nonlinear, are linearized using the Newton method or the $L$-scheme. 

In contrast, the alternate linearized splitting (AltLinS) schemes perform only one linearization step per iteration, see Figure \ref{scheme3}.
One iteration in the alternate splitting scheme can be written as:
\\Find  $\Psi_h^{n, j+1}, \theta_h^{n, j+1}$ such that
\begin{equation}
\begin{split}
\begin{cases} F^{lin}_1 (\Psi_h^{n, j+1}, \theta_h^{n,j+1}) &= 0,\\
 F^{lin}_2 (\Psi_h^{n, j+1}, \theta_h^{n, j+1},c_h^{n, j}) &= 0,
 \end{cases}\\
 \end{split}
\end{equation}
 and then $c_h^{n, j+1}$ such that
\begin{equation}
\begin{split}
F^{lin}_3 (\Psi_h^{n, j+1},\theta_h^{n, j+1}, c_h^{n, j+1}) = 0.
\end{split}
\end{equation}
Again, depending on which linearization is used, we refer to alternate splitting Newton (AltS-Newton) or alternate splitting $L$-scheme (AltS-LS). Both schemes will be presented in detail below.

In the following sections we will illustrate, in the details, the different schemes here investigated.

\subsubsection{The monolithic Newton method (Newton-Mono)}

The standard monolithic Newton method applied to \eqref{richards1} reads as:

\textbf{Problem P-Newton-Mono:} Let $ j>1 $ be fixed. Let $\Psi^{n-1}_h,\Psi^{n,j}_h, c_h^{n-1}, c_h^{n,j} \in V_h$, and $ \theta^{n-1}_h, \theta^{n,j}_h\in W_h$ be given, find $\Psi^{n,j+1}_h, c_h^{n,j+1} \in V_h$, and $\theta^{n,j+1}_h \in W_h$ such that
\begin{equation}\label{Newtonrichards}
\begin{split}
&<\theta^{n,j+1}_h - \theta^{n-1}_h, v_{1,h} > +\Delta t <K(\theta^{n,j}_h) (\nabla(\Psi^{n,j+1}_h)
+ \ez), \nabla v_{1,h} >  \\
+\Delta t &<\partial_\theta K(\theta^{n,j}_h) (\nabla(\Psi^{n,j}_h)+\ez)(\theta^{n,j+1}_h-\theta^{n,j}_h), \nabla v_{1,h} > = \Delta t <\mathbb{S}_1,v_{1,h} >
\end{split}
\end{equation}
\begin{equation}\label{psiNewton}
\begin{split}
\Delta t<\Psi^{n,j+1}_h , &w_{1,h}> = -\Delta t <p_{cap}(\theta^{n,j}_h,c^{n,j}_h), w_{1,h}>  - \Delta t <\partial_\theta p_{cap}(\theta^{n,j}_h,c^{n,j}_h) (\theta^{n,j+1}_h - \theta^{n,j}_h), w_{1,h}> \\
&-\Delta t <\partial_c p_{cap}(\theta^{n,j}_h,c^{n,j}_h) (c^{n,j+1}_h - c^{n,j}_h), w_{1,h}> 
+ <\tau(\theta^{n,j}_h) (\theta^{n,j+1}_h - \theta^{n-1}_h), w_{1,h}>\\
& +<\partial_\theta\tau(\theta^{n,j}_h) (\theta^{n,j}_h - \theta^{n-1}_h) (\theta^{n,j+1}_h - \theta^{n,j}_h), w_{1,h}>
+\Delta t \gamma <\Phi\left(\frac{\theta^{n,j}_h - \theta^{n-1}_h}{\Delta t}\right), w_{1,h}>\\
\end{split}
\end{equation}
and 
\begin{equation}\label{Newtontransp}
\begin{split}
< \theta^{n,j}_h &(c^{n,j+1}_h - c^{n-1}_h) + c^{n,j}_h (\theta^{n,j+1}_h - \theta^{n-1}_h), v_{2,h} > 
  + \Delta t < D\nabla c^{n,j+1}_h + \uwo c^{n,j+1}_h, \nabla v_{2,h}>  \\
  &+ \Delta t<R(c^{n,j}_h), v_{2,h}>  + \Delta t < \partial_c R(c^{n,j}_h)  (c^{n,j+1}_h - c^{n,j}_h)>  =\Delta t <\mathbb{S}_2,v_{2,h}>
\end{split}
\end{equation}
hold true for all $v_{1,h}, v_{2,h} \in V_h$, and for all $w_{1,h}\in W_h$. By $\partial_\theta$ we denote the partial derivative with respect to the water content $\theta$, and by $\partial_c$ the partial derivative with respect to the concentration $c$, and $\uwo:=-K(\theta^{n,j}_h)\nabla(\Psi^{n,j}_h+\ez)$. 

\subsubsection{The monolithic $L$-scheme (LS-Mono)}

The monolithic $L$-scheme for solving \eqref{richards1} reads:
\\ \textbf{Problem P-LS-Mono:} Let $ j>1 $ be fixed. Let $\Psi^{n-1}_h,\Psi^{n,j}_h, c^{n-1}_h, c_h^{n,j} \in V_h$, and $ \theta^{n-1}_h, \theta^{n,j}_h\in W_h$ be given, find $\Psi^{n,j+1}_h, c^{n,j+1}_h \in V_h$, and $\theta^{n,j+1}_h \in W_h$ such that

\begin{equation}\label{Lrichards}
<\theta^{n,j+1}_h - \theta^{n-1}_h, v_{1,h} >  +
\Delta t <K(\theta^{n,j}_h) (\nabla(\Psi^{n,j+1}_h)+\ez),\nabla v_{1,h} > +  L_1<\Psi^{n,j+1}_h-\Psi^{n,j}_h,v_{1,h}> = \Delta t <\mathbb{S}_1,v_{1,h} >
\end{equation}
\begin{equation}\label{Lpsi}
\begin{split}
\Delta t<\Psi^{n,j+1}_h , w_{1,h}> &= -\Delta t <p_{cap}(\theta^{n,j}_h,c^{n,j}_h), w_{1,h}> + <\tau(\theta^{n,j}_h) (\theta^{n,j+1}_h - \theta^{n-1}_h), w_{1,h}> \\
&+\Delta t \gamma <\Phi\left(\frac{\theta^{n,j}_h - \theta^{n-1}_h}{\Delta t}\right), w_{1,h}> + L_2 <(\theta^{n,j+1}_h - \theta^{n,j}_h) , w_{1,h}>\\
\end{split}
\end{equation}
and 
\begin{equation}\label{Ltransp}
\begin{split}
< \theta^{n,j}_h &(c^{n,j+1}_h - c^{n-1}_h) + c^{n,j}_h (\theta^{n,j+1}_h - \theta^{n-1}_h), v_{2,h} > 
 +\Delta t< D\nabla c^{n,j+1}_h + \uwo c^{n,j+1}_h, \nabla v_{2,h}>  \\
 &+\Delta t <R(c^{n,j}_h), v_{2,h}>  +  L_3<c^{n,j+1}_h-c^{n,j}_h,  v_{2,h}>
 =\Delta t<\mathbb{S}_2,v_{2,h}>
\end{split}
\end{equation}
hold true for all $v_{1,h}, v_{2,h} \in V_h$, and for all $w_{1,h}\in W_h$. $L_1,L_2$ and $L_3$ are three positive, user-defined parameters on which only mild conditions are imposed. We refer to \cite{Illiano2020,List2016,Pop2004} for the analysis of the numerical schemes which have inspired the ones presented here. Often, one needs to properly tune these parameters to obtain a robust and relatively fast solver.

\subsubsection{The nonlinear splitting approach (NonLinS)}

The nonlinear splitting approach for solving \eqref{richards1} reads:

\textbf{Problem P-NonLinS:} Let $ j>1 $ be fixed. Let $\Psi^{n-1}_h,\Psi^{n,j}_h, c^{n-1}_h, c^{n,j}_h \in V_h$ and $\theta^{n-1}_h, \theta^{n,j}_h\in W_h$ be given, find $\Psi^{n,j+1}_h \in V_h$, and $\theta^{n,j+1}_h \in W_h$ such that
\begin{equation}\label{NonLinrichards}
<\theta^{n,j+1}_h - \theta^{n-1}_h, v_{1,h} > +\Delta t <K(\theta^{n,j+1}_h) (\nabla(\Psi^{n,j+1}_h)+ \ez), \nabla v_{1,h} > =\Delta t <\mathbb{S}_1,v_{1,h} >
\end{equation}
\begin{equation}\label{NonLinpsi}
\begin{split}
\Delta t<\Psi^{n,j+1}_h , w_{1,h}> = -\Delta t <&p_{cap}(\theta^{n,j+1}_h,c^{n,j}_h), w_{1,h}> + <\tau(\theta^{n,j+1}_h) (\theta^{n,j+1}_h - \theta^{n-1}_h), w_{1,h}>\\
&+\Delta t \gamma <\Phi\left(\frac{\theta^{n,j}_h - \theta^{n-1}_h}{\Delta t}\right), w_{1,h}>
\end{split}
\end{equation}
holds true for all $v_{1,h} \in V_h$ and for all $w_{1,h}\in W_h$.

Then let $\Psi^{n-1}_h,\Psi^{n,j}_h, c^{n-1}_h, c^{n,j}_h \in V_h$ and $\theta^{n-1}_h, \theta^{n,j}_h\in W_h$ be given, $\Psi^{n,j+1}_h\in V_h$ and $\theta^{n,j+1}_h\in W_h$ are obtained from the equations above, find $c^{n,j+1}_h \in V_h$ such that
\begin{equation}\label{NonLintransp}
\begin{split}
< \theta^{n,j+1}_h (c^{n,j+1}_h - c^{n-1}_h) + c^{n,j}_h &(\theta^{n,j+1}_h - \theta^{n-1}_h), v_{2,h} >  
 + \Delta t < D\nabla c^{n,j+1}_h + \uwj c^{n,j+1}_h, \nabla v_{2,h}> \\
 &+\Delta t <R(c^{n,j+1}_h), v_{2,h}>  =\Delta t <\mathbb{S}_2,v_{2,h}>
\end{split}
\end{equation}
holds true for all $v_{2,h} \in V_h$. The water flux is given by $\uwj :=-K(\theta^{n,j+1}_h)\nabla(\Psi^{n,j+1}_h+\ez)$.

Observe that \eqref{NonLinrichards}--\eqref{NonLinpsi} and \eqref{NonLintransp} are nonlinear. To approximate their respective solutions, one can employ, e.g., the Newton method (NonLinS-Newton) or the $ L $-scheme (NonLinS-LS).

\subsubsection{The alternate splitting Newton method (Newton-AltLinS)}

Applied to \eqref{richards1}, the alternate splitting Newton method reads:

\textbf{Problem P-Newton-AltLinS:} Let $ j>1 $ be fixed. Let $\Psi^{n-1}_h,\Psi^{n,j}_h, c^{n-1}_h, c^{n,j}_h \in V_h$ and $\theta^{n-1}_h, \theta^{n,j}_h\in W_h$ be given, find $\Psi^{n,j+1}_h \in V_h$, and $\theta^{n,j+1}_h \in W_h$ such that
\begin{equation}\label{AltNewtrichards}
\begin{split}
&<\theta^{n,j+1}_h - \theta^{n-1}_h, v_{1,h} > 
+\Delta t <K(\theta^{n,j}_h) (\nabla(\Psi^{n,j+1}_h)
 + \ez), \nabla v_{1,h}  >  \\
+\Delta t &<\partial_\theta K(\theta^{n,j}_h) (\nabla(\Psi^{n,j}_h)
 +\ez)(\theta^{n,j+1}_h-\theta^{n,j}_h), \nabla v_{1,h}  >
 = \Delta t <\mathbb{S}_1,v_{1,h}  >
\end{split}
\end{equation}
\begin{equation}\label{AltNewtonpsi}
\begin{split}
\Delta t<\Psi^{n,j+1}_h , w_{1,h}> = -\Delta t <p_{cap}&(\theta^{n,j}_h,c^{n,j}_h), w_{1,h}>  - \Delta t <\partial_\theta p_{cap}(\theta^{n,j}_h,c^{n,j}_h) (\theta^{n,j+1}_h - \theta^{n,j}_h), w_{1,h}> \\
+ <\tau(\theta^{n,j}_h) (\theta^{n,j+1}_h - &\theta^{n-1}_h), w_{1,h}> 
+<\partial_\theta\tau(\theta^{n,j}_h) (\theta^{n,j}_h - \theta^{n-1}_h) (\theta^{n,j+1}_h - \theta^{n,j}_h), w_{1,h}>\\
&+\Delta t \gamma <\Phi\left(\frac{\theta^{n,j}_h - \theta^{n-1}_h}{\Delta t}\right), w_{1,h}>
\end{split}
\end{equation}
hold true for all $v_{1,h} \in V_h$ and $ w_{1,h} \in W_h$.

Then, with given $\Psi^{n-1}_h,\Psi^{n,j}_h, c^{n-1}_h, c^{n,j}_h \in V_h$ and $\theta^{n-1}_h, \theta^{n,j}_h\in W_h$, $\Psi^{n,j+1}_h\in V_h$ and $\theta^{n,j+1}_h\in W_h$ are obtained from the equations above, find $c^{n,j+1}_h \in V_h$ such that
\begin{equation}\label{AltNewttransp}
\begin{split}
< \theta^{n,j+1}_h &(c^{n,j+1}_h - c^{n-1}_h) + c^{n,j}_h (\theta^{n,j+1}_h - \theta^{n-1}_h), v_{2,h} >  + \Delta t < D\nabla c^{n,j+1}_h +\uwj c^{n,j+1}_h, \nabla v_{2,h}>  \\
  &+ \Delta t <R(c^{n,j}_h), v_{2,h}>  + \Delta t  <\partial_c R(c^{n,j}_h) (c^{n,j+1}_h - c^{n,j}_h), v_{2,h}>   = \Delta t <\mathbb{S}_2,v_{2,h}>
\end{split}
\end{equation}
hold true for all $v_{2,h} \in V_h$.

\subsubsection{The alternate splitting $L$-scheme (LS-AltLinS)}

The alternate splitting $L$-scheme for solving \eqref{richards1} is:

\textbf{Problem P-LS-AltLinS:} Let $ j>1 $ be fixed. Let $\Psi^{n-1}_h,\Psi^{n,j}_h, c^{n-1}_h, c^{n,j}_h \in V_h$ and $\theta_h^{n-1}, \theta^{n,j}_h\in W_h$ be given, find $\Psi^{n,j+1}_h \in V_h$, and $\theta^{n,j+1}_h \in W_h$ such that
\begin{equation}\label{LAltrichards}
\begin{split}
<\theta^{n,j+1}_h - \theta^{n-1}_h, v_{1,h} >  + \Delta t <K(\theta^{n,j}_h) (\nabla(\Psi^{n,j+1}_h)+\ez),\nabla v_{1,h} > +  L_1<\Psi^{n,j+1}_h-\Psi^{n,j}_h,v_{1,h} > = \Delta t <\mathbb{S}_1,v_{1,h}  >
\end{split}
\end{equation}
\begin{equation}\label{LAltpsi}
\begin{split}
\Delta t<\Psi^{n,j+1}_h , &w_{1,h}> = -\Delta t <p_{cap}(\theta^{n,j}_h,c^{n,j}_h), w_{1,h}> + <\tau(\theta^{n,j}_h) (\theta^{n,j+1}_h - \theta^{n-1}_h), w_{1,h}> \\
&+  L_2<\theta^{n,j+1}_h-\theta^{n,j}_h,w_{1,h} >+\Delta t \gamma <\Phi\left(\frac{\theta^{n,j}_h - \theta^{n-1}_h}{\Delta t}\right), w_{1,h}>
\end{split}
\end{equation}
holds true for all $v_{1,h}\in V_h$ and $w_{1,h} \in W_h$.

Then, with given $\Psi^{n-1}_h,\Psi^{n,j}_h, c^{n-1}_h, c^{n,j}_h \in V_h$ and $\theta^{n-1}_h, \theta^{n,j}_h\in W_h$, and $\Psi^{n,j+1}_h\in V_h$ and $\theta^{n,j+1}_h\in W_h$ from the equations above.  We find $c^{n,j+1}_h \in V_h$ such that
\begin{equation}\label{LAlttransp}
\begin{split}
< \theta^{n,j+1}_h &(c^{n,j+1}_h - c^{n-1}_h) + c^{n,j}_h (\theta^{n,j+1}_h - \theta^{n-1}_h), v_{2,h} >
+\Delta t< D\nabla c^{n,j+1}_h + \uwj c^{n,j+1}_h, \nabla v_{2,h}> \\
&+ \Delta t <R(c^{n,j}_h), v_{2,h}> + L_3 < c^{n,j+1}_h - c^{n,j}_h,  v_{2,h}>  =\Delta t<\mathbb{S}_2, v_{2,h}>
\end{split}
\end{equation}
hold true for all $v_{2,h} \in V_h$.

\begin{remark}
There exist multiple improved formulations of both the Newton method and L-scheme. We refer, among others, to the trust region techniques \cite{wang2013}, and the modified L-scheme in \cite{Mitra}.
\end{remark}

\begin{remark}(Stopping criterion)\label{remark2}
For all schemes (monolithic or splitting), the iterations are stopped when, 
\begin{equation*}
\norm{\Psi_h^{n,j+1}-\Psi_h^{n,j}}_\infty \leq \epsilon_1, \qquad \norm{\theta_h^{n,j+1}-\theta_h^{n,j}}_\infty \leq \epsilon_2 \qquad \text{and} \qquad \norm{c_h^{n,j+1}-c_h^{n,j}}_\infty \leq \epsilon_3,
\end{equation*} 
where by $\norm{\cdot}_\infty$ we mean the $ L^\infty (\Omega)$ norm.
Later on, for ease of presentation, we consider
 $\epsilon_1=\epsilon_2=\epsilon_3=\epsilon$. The parameter $\epsilon$ will be defined in the numerical section.
\end{remark}

\subsection{Anderson acceleration}\label{AndersonSec}

Although the L-scheme is robust and converges under mild restrictions, the convergence rate depends strongly on the linearization parameters. We refer to \cite{List2016,Pop2004,Slodicka} for the analysis in case of standard Richards equation, and to \cite{Karpinski} for the nonstandard model.
Tuning the parameters to obtain optimal results in terms of numbers of iterations and thus of computational times, can be tedious and time-consuming. The Anderson Acceleration (AA) is a powerful post-processing tool which can drastically reduce the numbers of iterations required by linearly convergent schemes, such as the L-scheme here investigated. In addiction, it reduces the need for finding close to optimal linearization parameters.

D. G. Anderson introduced the acceleration tool in 1965 \cite{Anderson1965}, and since then it has been investigated in multiple works, to name a few \cite{Both2019,Evans2019,Walker}.
We recall here the definition of AA, presented in \cite{Walker}, formulated for a general fixed point problem, of the form: given $g:\mathbb{R}^n \rightarrow \mathbb{R}^n$, solve the system $x = g(x)$.

\begin{algorithm}
\caption{Classical Fixed-Point iteration}
\begin{algorithmic}[1]
\State Given $x_0$
\For {$k=0,1,$... until convergence} 
\State $x_{k+1}=g(x_k)$
\EndFor
\end{algorithmic}
\end{algorithm}
Opposed to utilize only the last iteration $ x_k $, in the AA the new approximation is a linear combination of previously computed ones, see Algorithm 2.
 In the following, we denote by AA$(m)$ the Anderson acceleration where $m+1$ previously computed iterates are taken into account. With this, AA$(0)$ is the non-accelerated formulation.
\begin{algorithm}
\caption{Anderson Acceleration AA$(m)$}
\begin{algorithmic}[1]
\State Given $x_0$
\For {$k=1,2$... until convergence} 
\State Set $m_k = min\{m,k-1\}$
\State Define the matrix $F_k = (f_{k-m_k-1}, \cdots,f_{k-1})$, where $f_i=g(x_i)-x_i$
\State Find $\alpha  \in \mathbb{R}^{m_k+1}$ that solves
\begin{equation*}
\min\limits_{\alpha=(\alpha_0, \cdots, \alpha_{m_k})^T}\norm{F_k\alpha}\  \text{s.t.} \ \sum_{i=0}^{m_k}\alpha_i = 1.
\end{equation*}
\State Define $x_{k}:=\sum_{i=0}^{m_k} \alpha_{i}g(x_{k-m_k+i-1})$
\EndFor
\end{algorithmic}
\end{algorithm}
As revealed in the test cases below, this technique can drastically reduce the number of iterations required by the L-scheme.

The original formulation presented in \cite{Anderson1965} allows a for more general step, 
\begin{equation*}
x_{k}:=\beta_k\sum_{i=0}^{m_k} \alpha_{i}g(x_{k-m_k+i-1}) + (1-\beta_k) \sum_{i=0}^{m_k} \alpha_{i} x_{k-m_k+i-1},
\end{equation*}
for a user-defined tuning parameter $\beta_k\in(0,1]$. We considered the simplified formulation, obtained with $\beta_k=1$, because no improvements have been observed in the numerical results when using the extended one.

We remark that large values for the depth $m$ can result in an instability of the solution algorithm. When implementing the Anderson acceleration, one has to tune this parameter properly. A small $m$ could produce only a small reduction in the numbers of iterations; too large $m$ could result in a non-converging algorithm \cite{Fang2008}.

\begin{remark}
The definition of the nonlinear splitting solvers allows for different ways to apply the AA. We study three different loops: the coupling one and the linearizing ones, one for each set of equations. We apply the Anderson acceleration to each of them. Two different parameters, $m$ and $m_{lin}$, are defined. The former is used for the AA on the coupling loop, the latter for the implementation on the linearization ones. The same $m_{lin}$ will be used for the loop regarding the flow equations and for the one regarding the transport.
\end{remark}

\section{Numerical examples}
\label{numerical}

In the following, we consider four numerical examples with increasing complexity, based on a manufactured solution, and an example in which the boundary conditions drive the flow but no manufactured solution is given. The first four will differ in the different values for $\gamma$, $\delta$ and $\tau(\theta)$ taken into account. We have implemented the models and solving schemes in MRST, a toolbox based on Matlab for the simulations of flow in porous media \cite{mrst}. We use the two point flux approximation, one of the most common spatial discretization techniques. We remark that the linearization schemes and solving algorithms do not depend on the particular choice of the spatial discretization, so one may apply these solvers to other methods as well, without any difficulty.

The domain is the unit square $\Omega$ and the final time taken into consideration is $T=3$.
The simulations are performed on regular meshes, consisting of squares with sides $dx = 1/10,\ 1/20,$ and $1/40$.  The time steps are $\Delta t = T/25,\ T/50$ and $T/100$. The $L$ parameters, used in the L-scheme formulations, are $L_1=L_2=L_3 = 0.1$, if not specified otherwise. We took into consideration different values, but the aforementioned choice seems to produce a robust algorithm which required fewest iterations to achieve the convergence. For the ease of the presentation, we set the three $L$ parameters equal to each other; one could define different values for each parameter, investigating even further the linearization of each equation. We avoided this due to the application of the AA. We will observe that the schemes can be drastically accelerated, even though the $L$ parameters are not optimal.
 
The condition numbers, for the stiffness matrices resulting from the different solving algorithms are computed using the $L^1$ norm, and we here report the averaged values over the full simulation. A minus sign $(-)$, in the tables reporting on iterations and condition numbers, implies that the method failed to converge for the particular combination of the time step and mesh size. The tolerance $\epsilon$ used in the stopping criterion presented in Remark \ref{remark2} is $\epsilon=1e-6$. We always report the total numbers of iterations required by the full simulation, not the average number required by each time step. For the splitting solver, we present, separately, the condition numbers of both flow and transport equations. Furthermore, for the nonlinear splitting, the iterations are divided in two, the ones required by the flow equations and the ones for the transport. Finally, the condition numbers reported are obtained by averaging over the full simulations.

We apply the Anderson acceleration to each solving algorithm, always reporting the depths $m$ and $m_{lin}$ used. Once more, $m_{lin}$ is the Anderson parameter used for the acceleration of the linearization loops regarding the flow and transport equations in the nonlinear splitting solvers. 

Inspired by \cite{Lunowa2020}, the first four examples are constructed in such a way that the following is an exact solution:

\begin{equation}\label{theta_a}
\theta_m(x,y,t) = \begin{cases}  1-\frac{1}{2} \cos((t_1(x,y)-t)^2) &\mbox{if} \qquad t<t_1(x,y), \\
\frac{1}{2}&\mbox{if}\qquad t_1(x,y)\leq t\leq t_2(x,y),\\
1- \frac{1}{2} \cos((t-t_2(x,y))^2) \qquad &\mbox{if}\qquad t> t_2(x,y), 
\end{cases}
\end{equation}
\begin{equation}\label{psi_a}
\Psi_m(x,y,t) = \begin{cases} -p_{cap}(\theta_m) + \tau(\theta_m) \partial_t \theta_m - \gamma &\mbox{if}\qquad \partial_t \theta_m < -\delta,\\ 
-p_{cap}(\theta_m) + \tau(\theta_m) \partial_t \theta_m  + \frac{\gamma}{\delta} \partial_t \theta_m 
& \mbox{if}\qquad -\delta \leq \partial_t \theta_m \leq  \delta,\\
-p_{cap}(\theta_m) + \tau(\theta_m) \partial_t \theta_m + \gamma &\mbox{if}\qquad \partial_t \theta_m > \delta,
\end{cases}
\end{equation}
\begin{equation}\label{c_a}
c_m(x,y,t) = x(x-1) y(y-1) t,
\end{equation}
where $t_1(x,y) = xy$, $t_2(x,y) = xy+2$. Once the manufactured water content is defined, one obtains the pressure by simply using the second equation in (\ref{nonstandard_regularized}). The capillary pressure is expressed as $p_{cap}(\theta,c) = 1 - \theta^2 - 0.1\ c^3$ and the conductivity as $K(\theta) = 1 + \theta^2$. Even though such a formulation may appear non-realistic, we are mainly interested in the nonlinearities.  Furthermore, a nonlinear reaction term, $R(c) = c/(c+1)$, is taken into account in the transport equation and the diffusion/dispersion coefficient $D$ is set equal to $1$.

Given the analytical expressions above, we can easily define the initial conditions, the Dirichlet boundary conditions on the unit square and compute the source terms $\mathbb{S}_1$ and $\mathbb{S}_2$ such that $\Psi_m$, $\theta_m$ and $c_m$ are solutions of the system. 
In particular, the initial concentration and water content are
\begin{equation*}
\begin{alignedat}{2}
c(x,y,0) &= 0\qquad &\mbox{on}\qquad \Omega,\\
\theta(x,y,0) &=  1-\frac{1}{2} \cos\big(t_1(x,y)^2\big)\qquad &\mbox{on}\qquad \Omega.\\
\end{alignedat}
\end{equation*}
We impose a zero concentration $c(x,y,t) = 0$ on the boundary of the domain.
The remaining boundary conditions, concerning the pressure, are time dependent. One needs to compute $t_1$ and $t_2$ on every side of the unit square. Once the time intervals given by $t_1$ and $t_2$ are obtained, the pressure can easily be imposed.
For example, on the left side $x=0$, thus $t_1=0$ and $t_2=2$. The water content $\theta$ becomes
\begin{equation*}
\theta(0,y,t) = \theta_{left}(t) = \begin{cases} 
\frac{1}{2}&\mbox{if}\qquad 0< t< 2,\\
1- \frac{1}{2} \cos((t-2)^2) \qquad &\mbox{if}\qquad t\geq 2, 
\end{cases}
\end{equation*} 
and the resulting pressure boundary condition is
\begin{equation*}
\Psi_{left}(0,y,t) = \begin{cases} -p_{cap}(\theta_{left}) + \tau(\theta_{left}) \partial_t \theta_{left} - \gamma &\mbox{if}\qquad \partial_t \theta_{left} <-\delta,\\ 
-p_{cap}(\theta_{left}) + \tau(\theta_{left}) \partial_t \theta_{left}  + \frac{\gamma}{\delta} \partial_t \theta_{left} 
& \mbox{if}\qquad -\delta \leq \partial_t \theta_{left} \leq  \delta,\\
-p_{cap}(\theta_{left}) + \tau(\theta_{left}) \partial_t \theta_{left} + \gamma &\mbox{if}\qquad \partial_t \theta_{left} > \delta.
\end{cases}
\end{equation*}
Analogously, one can compute the pressure boundary conditions on the remaining sides.

\subsection{Example 1, $\gamma = 0$, $\tau(\theta)=1$}

In the first example we impose $\gamma = 0$, thus, the hysteresis effects are neglected but we include a dynamic effect by considering a constant $\tau(\theta)= 1$. 
We compare the different algorithms presented in Section \ref{section:algorithms}, reporting in the Tables \ref{tab:1.0dx} and \ref{tab:1.0dt} the total numbers of iterations required by each algorithm, and the condition numbers of the systems associated with each scheme. In the former, we investigate a fixed time step size, $\Delta t = T/25$, in the latter a fixed mesh size, $dx = 1/10$. 
As expected, a finer mesh results in worse conditioned systems, while smaller time steps give better conditioned ones.  Moreover, the total number of iterations is increasing as we reduce the time step; smaller $\Delta t$ implies more time steps and thus more iterations.

The schemes based on the L-scheme appear to be better conditioned than those based on the Newton method. The result is coherent with the theory \cite{Illiano2020,List2016,Mitra,Pop2004,Slodicka}. One could even improve the condition numbers by using larger $L$ parameters. However, larger values would have also increased the total numbers of iterations. 

The alternate splitting schemes are converging much faster than the nonlinear ones. It is also interesting to observe that the numbers of iterations, required by the alternate splitting schemes, are comparable with the ones associated with the monolithic solvers.  In \cite{Illiano2020}, we observed similar results when solving the models without hysteresis and dynamic effects.

We notice also some reduction in the number of iterations required by the L-schemes thanks to the Anderson acceleration. The results obtained for the non-accelerated L-schemes ($m=0$) are already optimal in terms of numbers of iterations; thus, the improvement can only be minimal. We report the total number of iterations for the full simulation, but mention that on average, for each time step, the L-scheme requires only five or six iterations. This is already a remarkable result, achieved thanks to the optimal L parameters. Furthermore, the Newton solvers have resulted in being slower when combined with the AA. This is coherent with the theory where it has been observed that quadratically convergent schemes, cannot be improved and the resulting accelerated solvers appear slower \cite{Evans2019}. 

\begin{table}[H]
\centering
\captionsetup{justification=centering,margin=2cm}
  \scalebox{.58}{
  \begin{tabular}{ c | c c | c c  c c  | c c c c}
\hline \hline
& Monolithic & &  & NonLinS &  & & & AltLinS &&\\
\hline
& Newton & & & Newton &  & & &Newton  & &\\
\hline
dx &\# iterations & condition \#  & \# iterations & cond. \# Flow && cond. \# Transport  & \# iterations & cond. \# Flow && cond. \# Transport\\
1/10  &  65  & 4.91e+02   & 56 - 50  &  1.82e+02   && 1.60e+02   &  66  &  1.83e+02  && 1.60e+02   \\ 
1/20  &  69  & 2.22e+03   & 57 - 50  &  7.61e+02   &&  6.33e+02  &  66  &  7.60e+02  &&  6.3359e+02  \\ 
1/40  &  70  & 1.11e+04   & 58 - 50  &   3.29e+03  &&  2.52e+03  &  66  & 3.42e+03   &&  2.5165e+03  \\
 \hline \hline
 & Newton &(AA m = 1) & & Newton (AA $m=m_{lin}$ = 1)&  & & &Newton (AA m = 1) & &\\
\hline
dx &\# iterations & condition \#  & \# iterations & cond. \# Flow && cond. \# Transport  & \# iterations & cond. \# Flow && cond. \# Transport\\
1/10  & 93   &  4.76e+02   & 65 - 50   & 1.82e+02  && 1.60e+02   &  74  &  1.82e+02  && 1.60e+02   \\ 
1/20  & 98   & 2.10e+03   &  66 - 50  &  7.54e+02  && 6.33e+02   &  74  &  7.59e+02  && 6.34e+02   \\ 
1/40  & 100   & 9.92e+03 & 69 - 50   & 3.28e+03   &&  2.52e+03   &  76  & 3.36e+03   &&  2.52e+03  \\
 \hline \hline
 & L-scheme & & & L-scheme &  & & & L-scheme & & \\
\hline
dx &\# iterations & condition \#  & \# iterations & cond. \# Flow && cond. \# Transport  & \# iterations & cond. \# Flow && cond. \# Transport\\ 
1/10  & 134  & 4.15e+02  & 117 - 116  & 1.43e+02  &&  1.36e+02  &  140  &  1.62e+02   &&  1.3717e+02  \\
1/20  & 140  & 1.82e+03  & 119 - 115  & 6.95e+02  &&  5.41e+02 & 144   &  7.37e+02   &&  5.41e+02  \\
1/40  & 146  &  8.52e+03 & 128 - 116  &  3.12e+03 && 2.15e+03   & 150   &  3.18e+03  &&   2.15 e03 \\
 \hline \hline
 & L-scheme & (AA m = 1)  & & L-scheme (AA $m=m_{lin}$ = 1) &  & & & L-scheme (AA m= 1) & & \\
\hline
dx &\# iterations & condition \#  & \# iterations & cond. \# Flow && cond. \# Transport  & \# iterations & cond. \# Flow && cond. \# Transport\\
1/10  &  127 &  4.14e+02 & 107 - 100  &  1.9690e+02    && 1.3681e+02   & 136   & 1.44e+02   && 1.33e+02   \\
1/20  &  129 &  1.85e+03 & 111 - 100  &  7.9408e+02    && 5.4122e+02   &  142  &  7.13e+02  && 5.31e+02   \\
1/40  &  130 &  8.80e+03 & 118 - 100  &   3.3661e+03   &&  2.1522e+03  &  146  &  3.01e+03  && 2.02e+03   \\
\hline\hline
\end{tabular}
}
\caption{Example 1: Total number of iterations and condition numbers for fixed $\Delta t=T/25$, and different $ dx $. Here, $L_1 = L_2 = L_3 = 0.1$ and $m=m_{lin}=1$.}
\label{tab:1.0dx}
\end{table}

\begin{table}[H]
\centering
\captionsetup{justification=centering,margin=2cm}
  \scalebox{.58}{
  \begin{tabular}{ c | c c | c c  c c  | c c c c}
\hline \hline
& Monolithic & &  & NonLinS &  & & & AltLinS &&\\
\hline
& Newton & & & Newton &  & & &Newton  & &\\
\hline
$\Delta t$&\# iterations & condition \#  & \# iterations & cond. \# Flow && cond. \# Transport  & \# iterations & cond. \# Flow && cond. \# Transport\\
T/25  & 65  &  4.91e+02 &56 - 50  &  1.82e+02   && 1.60e+02 & 66  &  1.83e+02  && 1.60e+02   \\  
T/50  & 103&  2.75e+02 & 98 - 50  &  1.23e+02   && 8.86e+01 & 99 & 1.24e+02   &&  8.89e+01  \\ 
T/100& 186& 1.97e+02  & 172 - 50  & 8.69e+01  && 4.73e+01 & 172& 8.74e+01 && 4.74e+01   \\
 \hline \hline
 & Newton &(AA m = 1) & & Newton (AA $m=m_{lin}$ = 1)&  & & &Newton (AA m = 1) & &\\
\hline
$\Delta t$&\# iterations & condition \#  & \# iterations & cond. \# Flow && cond. \# Transport  & \# iterations & cond. \# Flow && cond. \# Transport\\
T/25  & 93   & 4.76e+02  & 65 - 50   & 1.82e+02  && 1.60e+02 & 74  &  1.82e+02  && 1.60e+02   \\ 
T/50  & 137 & 2.75e+02  & 114 - 50  & 1.23e+02 &&8.86e+01   & 114&  1.23e+02  &&  8.89e+01  \\ 
T/100& 244 & 1.96e+02  & 201 - 100  & 8.59e+01 && 4.73e+01& 200   & 8.69e+01   && 4.74e+01   \\
 \hline \hline
 & L-scheme & & & L-scheme &  & & & L-scheme & & \\
\hline
$\Delta t$&\# iterations & condition \#  & \# iterations & cond. \# Flow && cond. \# Transport  & \# iterations & cond. \# Flow && cond. \# Transport\\ 
T/25   & 134  & 4.16e+02  & 117 - 116  & 1.43e+02  &&  1.36e+02 & 140  &  1.62e+02   &&  1.3717e+02  \\
T/50   & 219  & 2.39e+02  & 182 - 200   &  1.12e+02    && 7.51e+01&  218  &  1.2460e+02  &&   7.5203e+01 \\
T/100 & 425  & 1.64e+02 &  346 - 400   &  8.35e+01    && 3.98e+01&  438  &  8.2955e+01  &&  3.9923e+01  \\
 \hline \hline
 & L-scheme & (AA m = 1)  & & L-scheme (AA $m=m_{lin}$ = 1) &  & & & L-scheme (AA m = 1) & & \\
 \hline
$\Delta t$&\# iterations & condition \#  & \# iterations & cond. \# Flow && cond. \# Transport  & \# iterations & cond. \# Flow && cond. \# Transport\\
T/25  & 127  & 4.14e+02 & 107 - 100  &  1.97e+02    && 1.37e+02   &  136   & 1.44e+02   && 1.33e+02   \\
T/50  & 217  & 2.36e+02 & 192 - 160  &  1.25e+02    &&  7.51e+01  &  238  & 1.14e+02   &&   7.50e+01 \\
T/100 & 387  & 1.66e+02  & 347 - 300  &  8.33e+01    && 3.98e+01 & 432   & 8.21e+01   &&   3.91e+01 \\
\hline\hline
\end{tabular}
}
\caption{Example 1: Total number of iterations and condition numbers for fixed $dx=1/10$, and different $ \Delta t $. Here  $L_1 = L_2 = L_3 = 0.1$, $m=m_{lin}=1$.}
\label{tab:1.0dt}
\end{table}

In Table \ref{tab:1EOC} we present the numerical errors and the estimated order of convergence of the spatial discretization based on the successively refined meshes investigated. Given the manufactured solution $\Psi_m$, we compute the numerical error $e_\Psi=\norm{\Psi_m - \Psi_n}$, where $\Psi_n$ is the numerical pressure computed. Similarly, we can define $e_\theta$ and $e_c$. Furthermore, $e_{\Psi,1}$ is the numerical error obtained for the mesh size $dx=1/10$ and $\Delta t = T/25$, $e_{\Psi,2}$ for  $dx=1/20$ and $\Delta t = T/50$, and finally $e_{\Psi,3}$ for  $dx=1/40$ and $\Delta t = T/100$. $EOC=\log\left(\frac{e_{\Psi,i}}{e_{\Psi,i+1}}\right)/\log(2)$ is the estimated order of convergence. These results are independent from the solving algorithm taken into account, only the discretization approach plays a role. In this case we use a TPFA which is known to have a order of convergence equal to 1, as also reported here in the table.

\begin{table}[h!]
\begin{center}
	\captionsetup{justification=centering,margin=2cm}
\begin{tabular}{|l|l|l|l|l|l|}
\hline
&$e_{1}$ & EOC & $e_{2}$ & EOC & $e_{3}$ \\ \hline
$\Psi$     & 4.61e-02   &  0.98    &   2.33e-02    &   1.00   &    1.16e-02   \\ \hline
 $\theta$ &  1.31e-02   &   0.97    &  6.71e-03     &   0.99     &   3.38e-03       \\ \hline
$c$     &  6.24e-03  &  1.53    &    2.15e-03    &    1.32  &    8.60e-04     \\ \hline
\end{tabular}
\caption{Example 1: Numerical error and estimated order of convergence (EOC) of the discretization method.}
\label{tab:1EOC}
\end{center}
\end{table}

In Table \ref{tab:MsEx1} we tested different values of $m$ and L. We can observe that, for large L, the Monolithic L-Scheme, here investigated, requires more iterations than for smaller parameters. If many iterations are required to achieve the convergence at each time step, one can take in consideration larger $m$ values. For $L= 0.1$, the optimal choice, in terms of numbers of iterations, is $m=1$; larger values result in slightly slower schemes. For the largest $L$ tested, $L=2.3$, the optimal choice is $m=2$. Such L value corresponds to the theoretical L, $\max\limits_{\theta}\{\frac{\partial p_{cap}}{\partial \theta}\}\approx2.3$. To ensure the monotone convergence of the scheme it has been proved, that the parameter chosen must be larger than the Lipschitz constant of the nonlinearity, in this case the capillary pressure as a function of $ \theta $ (see \cite{List2016,Pop2004,Slodicka}). We set $L_1$, $ L_2 $ and $L_3$ equal to the theoretical L computed for the capillary pressure.

We can conclude that it is possible to obtain significant improvements by investigating the AA and thus finding the appropriate depth $m$. In this work, we have focused more on individuating the optimal L parameters, as refining the AA can be done more easily. The depths used are, in fact, small natural numbers.

\begin{table}[h!]
\begin{center}
\captionsetup{justification=centering,margin=2cm}
\begin{tabular}{|l|l|l|l|l|l|}
\hline
L &$m=0$& $m=1$ & $m=2$ & $m=3$ & $m=5$  \\\hline
.1   &	146  &	 130  & 146   & 179    &	 215 	\\\hline  
.5   &	 247 &	161   & 152	 &	  172   &	  200 	\\\hline   
1    &	411  &	202   & 168 &	 165   &  195 	\\\hline 
2    &	711  &	290   & 199 	 &  206  &	 217 	\\\hline 
2.3 &	810 &	317   &  207	 & 215   &	227  	\\\hline 
\end{tabular}
\caption{Example 1: Comparison of number of iterations for different $m$ and L parameters for L-Mono. \\Here, $dx=1/40$, $\Delta t = T/25$ and $L_1=L_2=L_3=L$.}
\label{tab:MsEx1}
\end{center}
\end{table}

\begin{figure}[h!]
	\begin{center} 
		\captionsetup{justification=centering,margin=2cm}
	  \begin{subfigure}{.49\textwidth}
			\includegraphics[width=1\linewidth]{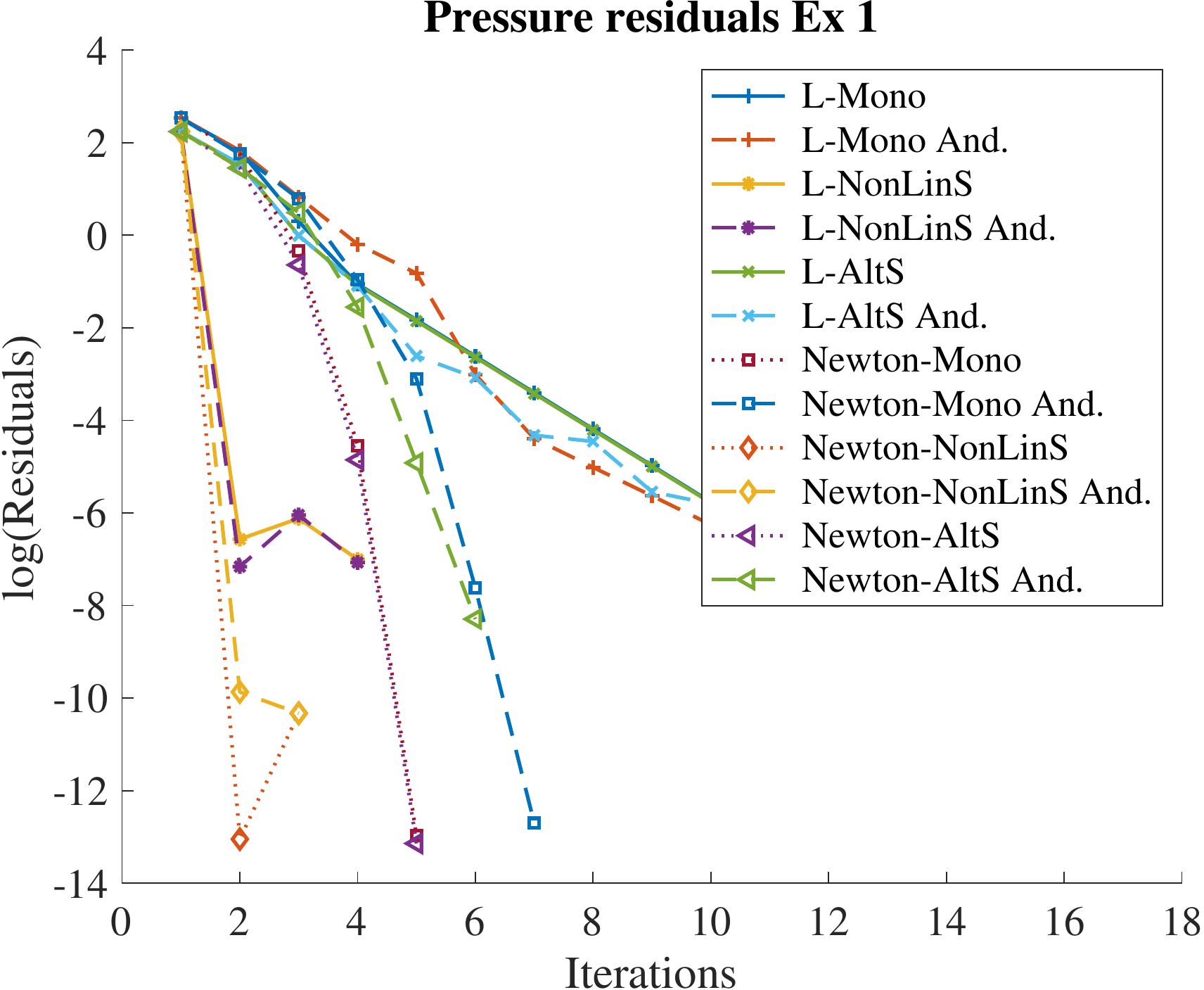}
			\caption{Pressure residuals.}
			\end{subfigure}
  \begin{subfigure}{.49\textwidth}
		  \includegraphics[width=1\linewidth]{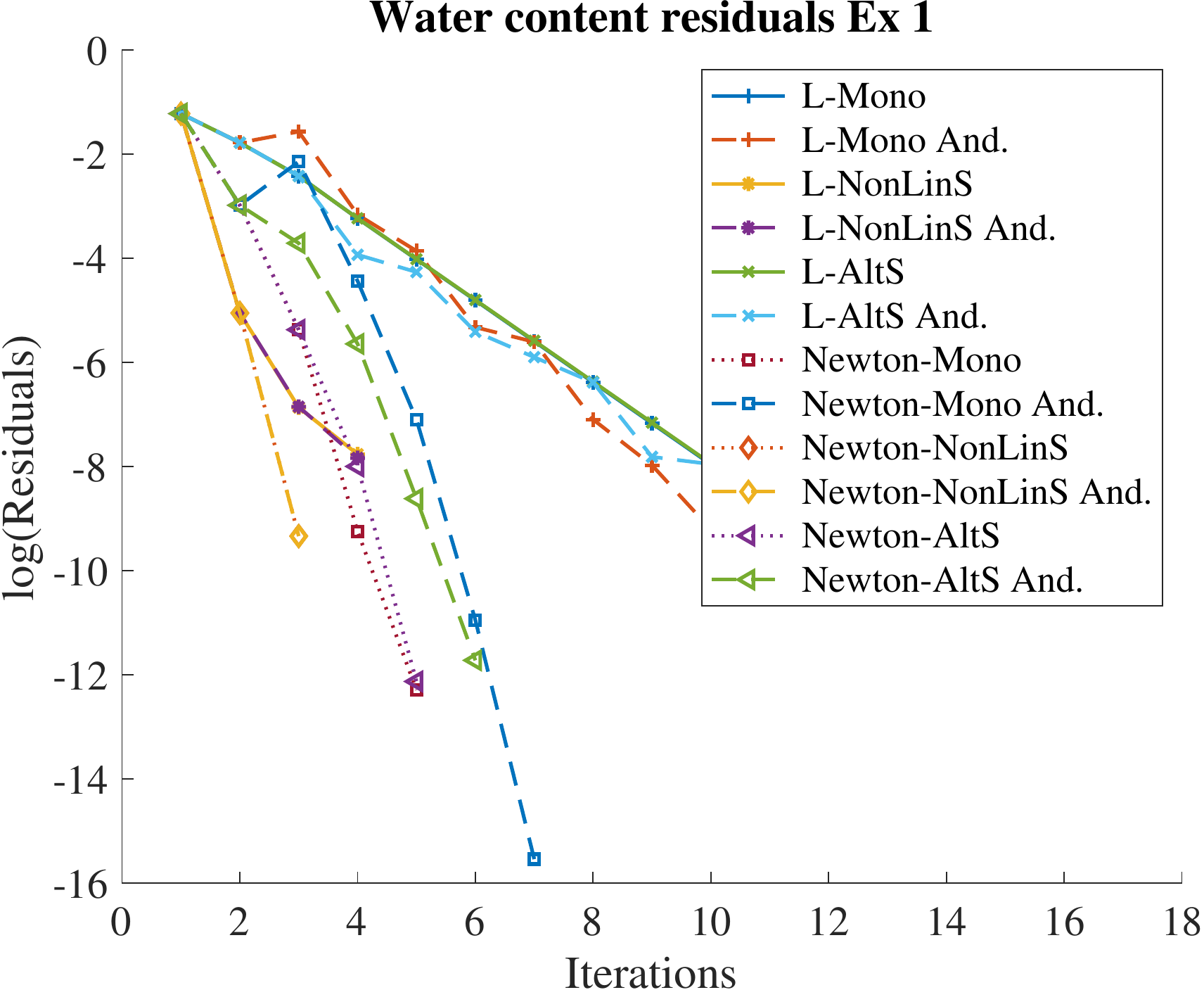}
			\caption{Water content residuals.}
			\end{subfigure}\\
			  \begin{subfigure}{.49\textwidth}
\includegraphics[width=1\linewidth]{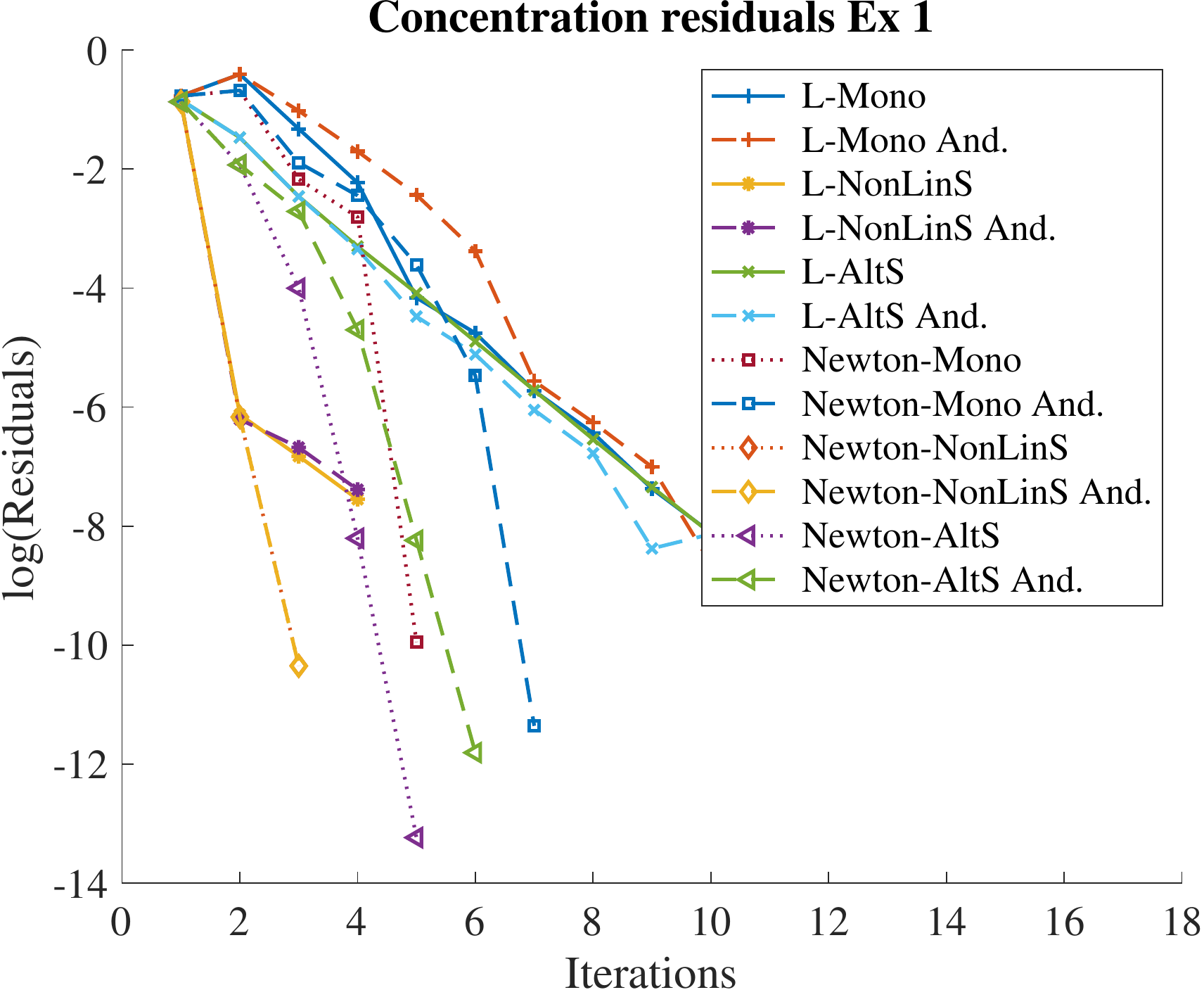}
			\caption{Concentration residuals.}
			\end{subfigure}		
		\caption{Example 1: Residuals of each unknown at the final time step, for the different schemes. \\Here, $L_1=L_2=L_3=0.1$, $m=m_{lin}=1$, $dx=1/40$, $\Delta t=T/25$.}
		\label{fig:Ex1Residuals}
	\end{center}
\end{figure}

Finally, we investigate the order of convergence of the linearization schemes. In Figure \ref{fig:Ex1Residuals}, we plot the residuals of pressure, water content and concentration, obtained at the final time step for the finest mesh size, $dx= 1/40$ and $\Delta t = T/25$. We can deduce the rates of convergence of the different linearization schemes. The L-schemes appear to be, as expected, linearly convergent in term of numbers of iterations. The AA improves the results only slightly, as already observed in Tables \ref{tab:1.0dx} and \ref{tab:1.0dt}. This is justified by the fact that the L parameters ($L_1,L_2,L_3$) chosen here, appear to be optimal. The Newton methods are instead quadratically convergent. More precise results are observable in Table \ref{tab:1.1EOC}. Here we present the exact order of convergence of the different schemes. Given the residual of each unknown ($res_\Psi,res_\theta,res_c$), at each iteration $j$, we compute the order of convergence as follow: $ORD_j=\big(\log(res_{j+1}/res_j)\big)/\big(\log(res_j/res_{j-1})\big)$. For a fixed time step, we can average the orders obtained over the number of iterations required to achieve the convergence. We report below the values obtained by investigating the final time step; similar results have been observed for previous time steps.

\begin{table}[h!]
\begin{center}
\begin{tabular}{|c|c|c|c|c|c|c|}
\hline
&\scriptsize{LS-Mono} & \scriptsize{LS-Mono And.} & \scriptsize{LS-NonLinS} & \scriptsize{LS-NonLinS And.}  & \scriptsize{LS-AltLinS} & \scriptsize{LS-AltLinS And.}\\ \hline
$\Psi$   &  1.00 & 1.26   &1.13    & 1.46  & 1.13 &2.07       \\ \hline
$c$        &  0.96 & 1.55  & 1.01   & 1.77   & 1.01 & 1.30   \\ \hline
$\theta$ & 1.01 & 1.40  & 0.89   &  1.01   & 0.94 &1.32  \\ \hline
 & \scriptsize{New.-Mono} & \scriptsize{New.-Mono And.}  & \scriptsize{New.-NonLinS} & \scriptsize{New.-NonLinS And.} & \scriptsize{New.-AltLinS} & \scriptsize{New.-AltLinS And.} \\   \hline
 $\Psi$   & 2.03  & 1.61 & 1.94   &  1.86  & 2.04 & 2.03\\ \hline
 $c$        &  2.41 &  2.01 & 1.61   &  1.41  &1.85 & 1.46\\ \hline
 $\theta$&  1.97 &  0.57 & 1.81   & 0.97    &1.61 & 0.95\\ \hline
\end{tabular}
\caption{Example 1: Estimated order of convergence for the different linearization schemes.}
\label{tab:1.1EOC}
\end{center}
\end{table}
 
\begin{remark}
Since the Anderson acceleration with small depth is a cheap post-processing step, reducing the number of iterations directly reduces the CPU time almost proportionally.
\end{remark}

\newpage
\subsection{Example 2, $\gamma=0$, $\tau(\theta)=1+ \theta^2$}

In the second example, the setup of the first is extended by adopting a nonlinear $\tau$, precisely $\tau(\theta) = 1+ \theta^2$. In Tables \ref{tab:2.0dx} and \ref{tab:2.0dt}, we present the condition numbers and the required iteration counts associated with each solving algorithm. 

The introduction of a nonlinear $\tau(\theta)$ increases the numbers of iterations required by each solver. The L-scheme is linearly convergent while the Newton method is quadratically convergent. Furthermore, the conclusions from Example 1 concerning the AA remain the same. In particular, we observe only small reductions in the numbers of iterations required by the L-schemes. Once more, this is justified by the optimal choice of the L parameters. We can observe that, for each time step, only a few iterations are required; thus, no further acceleration is expected.

\begin{table}[H]
\centering
\captionsetup{justification=centering,margin=2cm}
  \scalebox{.58}{
  \begin{tabular}{ c | c c | c c  c c  | c c c c}
\hline \hline
& Monolithic & &  & NonLinS &  & & & AltLinS &&\\
\hline
& Newton & & & Newton &  & & &Newton  & &\\
\hline
dx &\# iterations & condition \#  & \# iterations & cond. \# Flow && cond. \# Transport  & \# iterations & cond. \# Flow && cond. \# Transport\\
1/10  & 67  & 481.17  &  59 - 50 & 161.39 && 161.00  &  64     &    152.62 &&  161.76    \\
1/20  & 69  & 2.20e+03  & 57 - 50  & 656.26    && 636.33  &  66    &  644.98   &&  644.99    \\
1/40  & 70  & 1.10e+04  & 59 - 50  & 2.85e+03 && 2.53e+03  &    68  &  2.93e+03   &&  2.53e+03    \\
 \hline \hline
 & Newton & (Anderson m = 1)& & Newton (AA $m=m_{lin}$ = 1) &  & & &Newton (Anderson m = 1) & &\\
\hline
dx &\# iterations & condition \#  & \# iterations & cond. \# Flow && cond. \# Transport  & \# iterations & cond. \# Flow && cond. \# Transport\\
1/10  & 94 & 467.04 & 69 - 50  &  160.36   && 161.00  & 76  &   152.45 &&  161.94   \\
1/20  & 98 & 2.07e+03  & 68 - 50  &  636.14   && 636.33  & 76    &   642.25 &&   638.40  \\
1/40  & 100 & 9.85e+03  & 71 - 50  &  2.80e+03   &&  2.53e+03 &    80 &  2.87e+03  &&   2.53e+03  \\
 \hline \hline
 & L-scheme & & & L-scheme &  & & & L-scheme & & \\
 \hline
dx &\# iterations & condition \#  & \# iterations & cond. \# Flow && cond. \# Transport  & \# iterations & cond. \# Flow && cond. \# Transport\\ 
1/10  & 136 & 403.49  & 112 - 117  &  158.33   && 137.81  &  130   &  158.05  &&  138.50   \\
1/20  & 139 & 1.77e+03  & 120 - 116  &  644.85   && 544.34 &  136   &  651.07  &&  545.85   \\
1/40  & 144 & 8.34e+03  & 125 - 116  &  2.73e+03 && 2.16e+03 &142     & 2.81e+03   &&  2.16e+03   \\
 \hline \hline
 & L-scheme & (Anderson m = 1)  & & L-scheme (AA $m=m_{lin}$ = 1) &  & & & L-scheme (AA m = 1) & & \\
 \hline
dx &\# iterations & condition \#  & \# iterations & cond. \# Flow && cond. \# Transport  & \# iterations & cond. \# Flow && cond. \# Transport\\
1/10  & 129 & 401.61  & 106 - 100  &  158.87   &&  137.74 & 132    &  156.89  &&   138.21  \\
1/20  & 131 & 1.78e+03  & 113 - 100  & 644.63  && 543.86  &  134   &  647.35  &&  545.18   \\
1/40  & 137 & 8.44e+03  & 117 - 100  &  2.76e+03   && 2.16e+03  &  142   &  2.78e+03  && 2.16e+03    \\
\hline
\end{tabular}
}
\caption{Example 2: Total number of iterations and condition numbers for fixed $\Delta t=T/25$, and different $ dx $. Here, $L_1 = L_2 = L_3 = 0.1$ and $m=m_{lin}=1$.}
\label{tab:2.0dx}
\end{table}

\begin{table}[H]
\centering
\captionsetup{justification=centering,margin=2cm}
  \scalebox{.58}{
  \begin{tabular}{ c | c c | c c  c c  | c c c c}
\hline \hline
& Monolithic & &  & NonLinS &  & & & AltLinS &&\\
\hline
& Newton & & & Newton &  & & &Newton  & &\\
\hline
$\Delta t$&\# iterations & condition \#  & \# iterations & cond. \# Flow && cond. \# Transport  & \# iterations & cond. \# Flow && cond. \# Transport\\
T/25   & 67  & 481.17 & 59 - 50 & 161.39 && 161.00  &  64     &    152.62 &&  161.76    \\
T/50   & 107 & 256.90  & 102 - 50  & 95.40&&  89.26 &   102    &    95.08 &&   89.80   \\
T/100 & 193 & 159.33  & 184 - 100 &73.42&& 47.70  &   184    &    72.37 &&  47.91    \\
 \hline \hline
& Newton&  (AA m = 1)& & Newton (AA $m=m_{lin}$ = 1) &  & & &Newton (AA m = 1) & &\\
\hline
$\Delta t$&\# iterations & condition \#  & \# iterations & cond. \# Flow && cond. \# Transport  & \# iterations & cond. \# Flow && cond. \# Transport\\
T/25   & 94 & 467.04    & 69 - 50  &  160.36   && 161.00 & 76  &   152.45 &&  161.94   \\
T/50   & 142 & 254.28  & 117 - 50  &  94.64   &&  89.26 &  118   &   95.71 &&  89.88   \\
T/100 & 259 & 159.96  & 207 - 100  &  73.05   && 47.70  & 208   &   72.78 &&  47.95   \\
 \hline \hline
 & L-scheme & & & L-scheme &  & & & L-scheme & & \\
\hline
$\Delta t$&\# iterations & condition \#  & \# iterations & cond. \# Flow &&  cond. \# Transport  & \# iterations & cond. \# Flow && cond. \# Transport\\ 
T/25   & 136 & 403.49  & 112 - 117  &  157.86&&  137.81 & 130   &  158.05  &&  138.50   \\
T/50   & 220 & 211.60  & 185 - 201   & 98.63 &&  75.62  & 220    &   97.80 &&  75.86   \\
T/100 & 433 & 125.49  & 356 - 400   &  65.18 && 40.16   &  446   &   65.24 &&  40.28   \\
 \hline \hline
 & L-scheme & (AA m = 1)  & & L-scheme (AA $m=m_{lin}$ = 1) &  & & & L-scheme (AA m = 1) & & \\
 \hline
$\Delta t$&\# iterations & condition \#  & \# iterations & cond. \# Flow && cond. \# Transport  & \# iterations & cond. \# Flow && cond. \# Transport\\
T/25   & 129 & 401.61  & 106 - 100  &  158.87   &&  137.74 & 132    &  156.89  &&   138.21  \\
T/50   & 214 & 212.13  & 192 - 161  &  97.72   &&  75.69 &  240   &   97.63 &&  75.66   \\
T/100 & 388 & 126.15  &  355 - 302 &  64.73   && 40.16  &  440   &   64.45 && 40.17    \\
\hline
\end{tabular}
}
\caption{Example 2: Total number of iterations and condition numbers for fixed $dx=1/10$, and different $ \Delta t $. Here  $L_1 = L_2 = L_3 = 0.1$, $m=m_{lin}=1$.}
\label{tab:2.0dt}
\end{table}

As for the results presented in Tables \ref{tab:2.0dx} and \ref{tab:2.0dt}, the numerical errors and EOC reported in Table \ref{tab:2EOC} are similar to the ones from the first example. 

\begin{table}[h!]
\begin{center}
	\captionsetup{justification=centering,margin=2cm}
\begin{tabular}{|l|l|l|l|l|l|}
\hline
  &$e_{1}$ & EOC & $e_{2}$ & EOC & $e_{3}$ \\ \hline
 $\Psi$     &  0.0759    &  0.9688    &  0.0388    &  0.9913    &   0.0195      \\ \hline
 $\theta$  &  0.0140    & 0.9556  &    0.0072 &    0.9830 &    0.0036   \\ \hline
$c$       &  0.0084   &  1.5198    &    0.0029   &   1.3037    &     0.0012   \\ \hline
\end{tabular}
\caption{Example 2: Numerical error and estimated order of convergence (EOC) of the discretization method.}
\label{tab:2EOC}
\end{center}
\end{table}

\begin{figure}[h!]
	\begin{center} 
		\captionsetup{justification=centering,margin=2cm}
	  \begin{subfigure}{.49\textwidth}
			\includegraphics[width=1\linewidth]{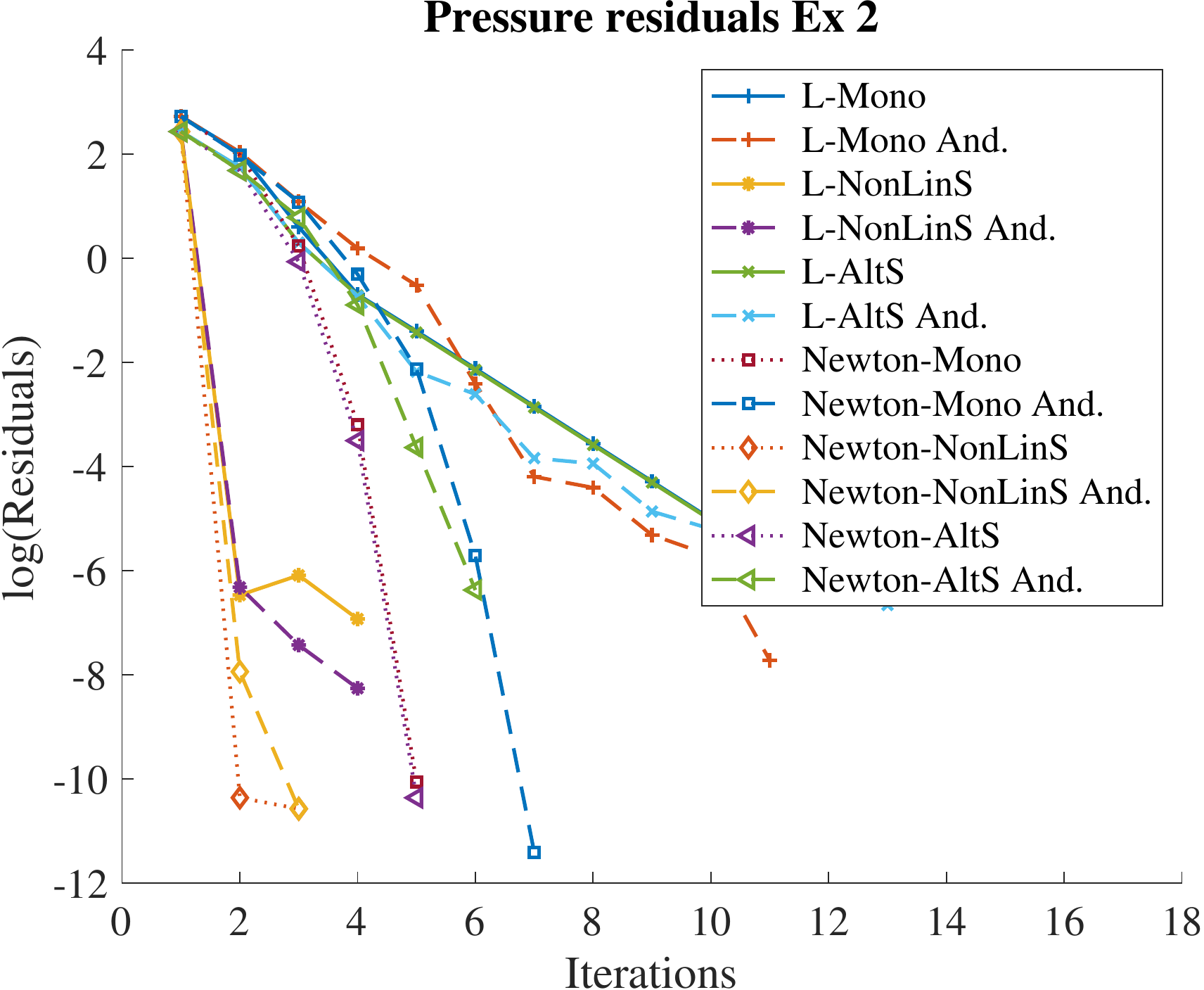}
			\caption{Pressure residuals.}
			\end{subfigure}
  \begin{subfigure}{.49\textwidth}
		  \includegraphics[width=1\linewidth]{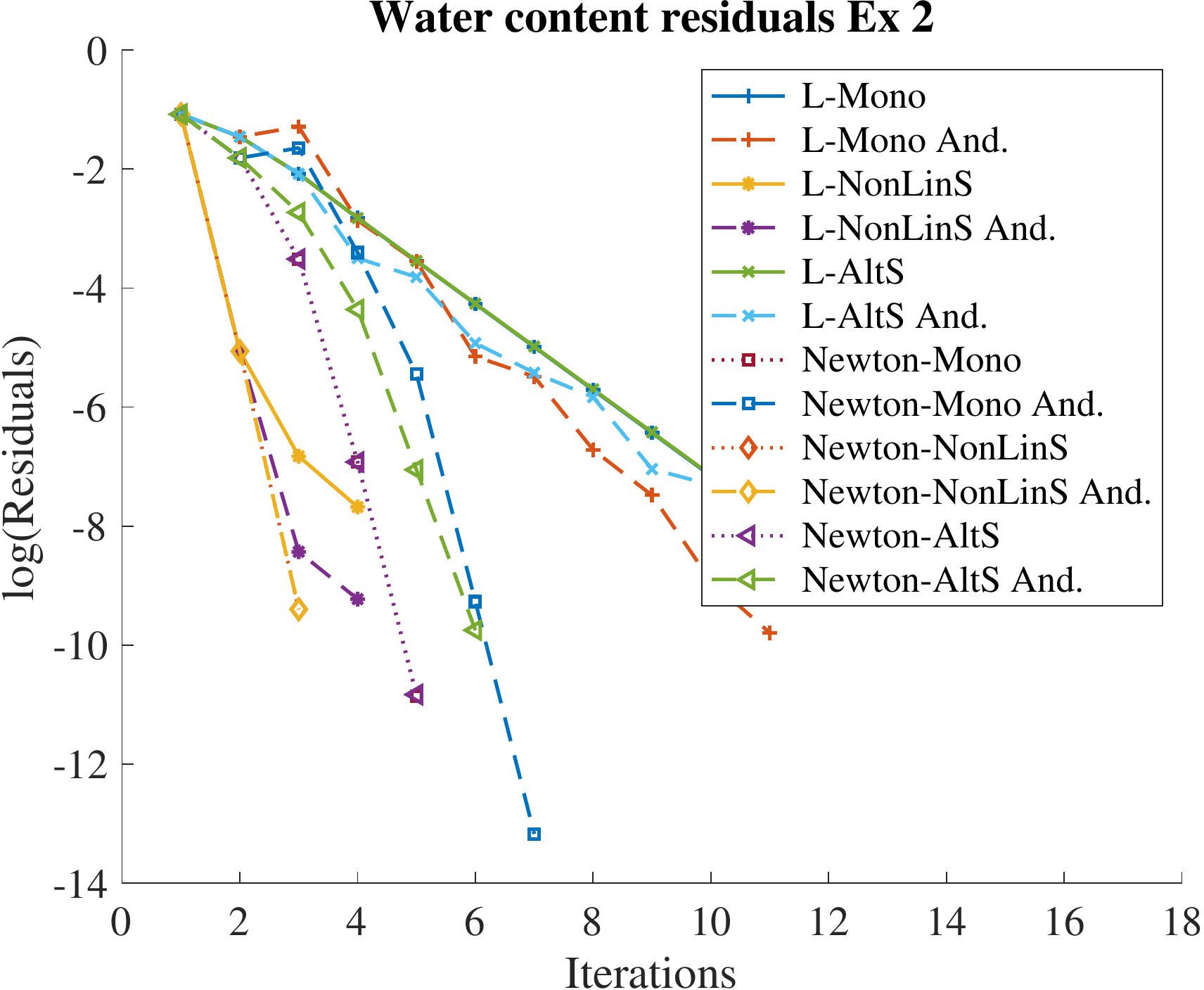}
			\caption{Water content residuals.}
			\end{subfigure}\\
			  \begin{subfigure}{.49\textwidth}
\includegraphics[width=1\linewidth]{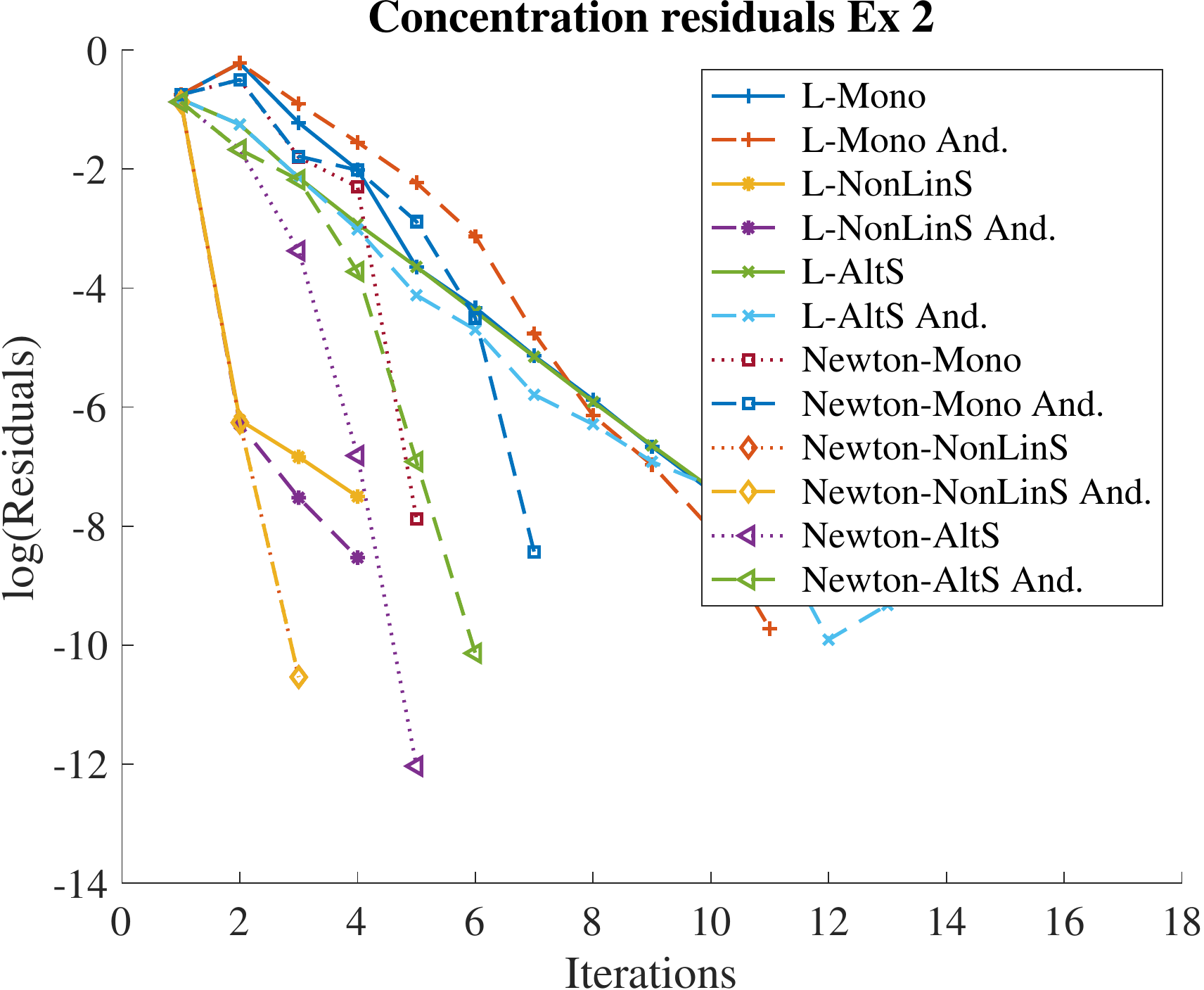}
			\caption{Concentration residuals.}
			\end{subfigure}		
		\caption{Example 2: Residuals of each unknown at the final time step, for the different schemes. Here, $L_1=L_2=L_3=0.1$, $m=m_{lin}=1$, $dx=1/40$, $\Delta t=T/25$. }
		\label{fig:Ex2Residuals}
	\end{center}
\end{figure}

In Figure \ref{fig:Ex2Residuals}, we report the residuals of the pressure, water content and concentration at the final time step. The L-schemes are linearly convergent, and applying the AA does not result in significant improvements. The convergence rates and number of iterations remain the same. Also for the Newton solvers, since they are quadratically convergent, the AA cannot improve this aspect. Table \ref{tab:2.1EOC} presents the precise order of convergence of the different linearization schemes.

\begin{table}[h!]
\begin{center}
	\captionsetup{justification=centering,margin=2cm}
\begin{tabular}{|c|c|c|c|c|c|c|}
\hline
&\scriptsize{LS-Mono} & \scriptsize{LS-Mono And.} & \scriptsize{LS-NonLinS} & \scriptsize{LS-NonLinS And.} & \scriptsize{LS-AltLinS} & \scriptsize{LS-AltLinS And.}\\ \hline
$\Psi$   &  1.07  &  1.40  &  1.11   &  1.36   &  1.14  & 1.24        \\ \hline
$c$        &  0.99  & 1.23  &  0.98   &  1.45    &  0.96 & 1.10      \\ \hline
$\theta$ & 1.03  &  1.15  &  0.97   &   1.25   &  0.93  & 0.98      \\ \hline
 & \scriptsize{New.-Mono} & \scriptsize{New.-Mono And.}  & \scriptsize{New.-NonLinS} & \scriptsize{New.-NonLinS And.} & \scriptsize{New.-AltLinS} & \scriptsize{New.-AltLinS And.} \\   \hline
 $\Psi$   &  1.61  & 1.58 &  1.97  &  1.69      &  2.15  & 2.14       \\ \hline
 $c$       & 2.68  &  1.47 &  1.98  &  1.54     &   1.98  & 1.83       \\ \hline
 $\theta$& 1.99 & 1.61  & 1.89  &  1.70     &   2.17  &  1.95      \\ \hline
\end{tabular}
\caption{Example 2: Estimated order of convergence for the different linearization schemes.}
\label{tab:2.1EOC}
\end{center}
\end{table}

\newpage
\subsection{Example 3, $\gamma=1$, $\delta=5e-3$, $\tau(\theta)=0$}
With the same manufactured solutions, we now consider the case without dynamic effects ($\tau(\theta)=0$), but include hysteresis by choosing $\gamma=1$ and $\delta=5e-3$. 

From the results in Tables \ref{tab:3.0dx} and \ref{tab:3.0dt}, we notice that the Newton method, in all its formulations, fails to converge. In Table \ref{tab:3.0dt}, smaller time steps are taken, but no improvements are observable. A further reduction of the time step could have resulted in converging Newton solvers but the total numbers of iterations for the full simulation would have been larger than the ones required by the L-schemes on fewer but larger time steps. In contrast, the L-schemes are more robust and, even though requiring a higher number of iterations than previously, they converge. We take $L_1=L_2=L_3=L=1$, which appear to be the optimal choice in terms of numbers of iterations. 

The AA improves the convergence of the L-schemes. This is the first example of this study in which the results obtained thanks to the AA are improved substantially. This is due to the presence of the hysteresis, requiring a large $L$ for the overall convergence, and thus the total numbers of iterations is larger. On average, the monolithic L-scheme solver requires circa 18 iterations per time step. For $m=1$, the AA reduces the iterations by circa $50\%$. Different $m$ values have been tested but none of the ones investigated lead to the convergence of the Newton schemes. On all tests, Newton has failed to converge, whereas the L-schemes converged and the AA yields further improvement.

\begin{table}[H]
\centering
\captionsetup{justification=centering,margin=2cm}
    \scalebox{.55}{
  \begin{tabular}{ c | c c | c c  c c  | c c c c}
\hline \hline
& Monolithic & &  & NonLinS &  & & & AltLinS &&\\
\hline
& Newton & & & Newton &  & & &Newton  & &\\
\hline
dx &\# iterations & condition \#  & \# iterations & cond. \# Flow && cond. \# Transport  & \# iterations & cond. \# Flow && cond. \# Transport\\
1/10  & - & - & - & - && -& - & - && -\\ 
1/20  & - & - & - & - && -& - & - && - \\ 
1/40  & - & - & - & - && -& - & - && - \\  
 \hline \hline
 & L-scheme & & & L-scheme &  & & & L-scheme & & \\
 \hline
dx &\# iterations & condition \#  & \# iterations & cond. \# Flow && cond. \# Transport  & \# iterations & cond. \# Flow && cond. \# Transport\\
1/10  & 448  & 409.16  &  210 - 441 &  484.83 && 69.22  &  450 &  361.76 && 69.34 \\ 
1/20  & 456  & 1.65e+03  & 266 - 439  & 1.98e+03  && 259.35  &  452 & 1.46e+03  &&  260.02\\ 
1/40  & 468  & 6.62e+03  &  276 - 438 &  7.74e+03 && 996.36  &  460 &  5.88e+03 &&  999.13\\ 
 \hline \hline
 & L-scheme & (AA m = 2)  & & L-scheme (AA m = 2, $m_{lin}$ = 5) &  & & & L-scheme (AA m = 1) & & \\
\hline
dx &\# iterations & condition \#  & \# iterations & cond. \# Flow && cond. \# Transport  & \# iterations & cond. \# Flow && cond. \# Transport\\
1/10  & 226  & 468.28  & 179 - 150  & 497.56  &&  70.10 &  328 &  450.28 && 71.81 \\ 
1/20  & 278 &  1.97e+03 &  187 - 141 & 2.03e+03  && 261.40  & 408   & 2.09e+03  && 269.72 \\ 
1/40  & 303  & 8.24e+03  & -  & -  && -  & 378  & 7.76e+03  && 967.29 \\
\hline
\end{tabular}
}
\caption{Example 3: Total number of iterations and condition numbers for fixed $\Delta t=T/25$, and different $ dx $. Here, $L_1 = L_2 = L_3 = 1$, different $m$ and $m_{lin}$ are taken into account.}
\label{tab:3.0dx}
\end{table}

\begin{table}[H]
\centering
\captionsetup{justification=centering,margin=2cm}
    \scalebox{.5}{
  \begin{tabular}{ c | c c | c c  c c  | c c c c}
\hline \hline
& Monolithic & &  & NonLinS &  & & & AltLinS &&\\
\hline
& Newton & & & Newton &  & & &Newton  & &\\
\hline
dx &\# iterations & condition \#  & \# iterations & cond. \# Flow && cond. \# Transport  & \# iterations & cond. \# Flow && cond. \# Transport\\
T/25   & - & - & - & -  && - & - & - && -\\ 
T/50   & - & - & - & -  && - & - & - && -\\ 
T/100 & - & - & - & -  && - & - & - && -\\
 \hline \hline
 & L-scheme & & & L-scheme &  & & & L-scheme & & \\
 \hline
&  & & &  &  cond. \# & & &  & cond. \# &\\
\hline
$\Delta t$ &\# iterations & condition \#  & \# iterations & Flow && Transport  & \# iterations & Flow && Transport\\ 
T/25   & 448  & 409.16   & 210 - 441 &  484.83 && 69.22 & 450 &  361.76 && 69.34 \\ 
T/50   & 836  &  363.15  & 513 - 846 &  422.16  && 35.96  & 838   &   332.02 && 36.03  \\ 
T/100 & 1787  & 395.55 & 1261 - 1597& 428.24&& 18.91  & 1764   &  363.40  && 19.07  \\
 \hline \hline
 & L-scheme & (AA m = 1)  & & L-scheme (AA m = 2, $m_{lin}$ = 5) &  & & & L-scheme (AA m = 1) & & \\
 \hline
dx &\# iterations & condition \#  & \# iterations & cond. \# Flow && cond. \# Transport  & \# iterations & cond. \# Flow && cond. \# Transport\\
T/25   & 226  & 468.28 & 179 - 150  & 497.56  &&  70.10  & 328 &  450.28 && 71.81 \\ 
T/50   & 533  & 410.93  & 346 - 316   & 504.00   && 36.85   &  664  &  442.54  &&  37.43 \\ 
T/100 & 1217& 467.07   & 861 - 944 ($m_{lin}$ = 1)  &  491.24  && 19.15  &  1842  & 473.32   && 19.75  \\ 
\hline
\end{tabular}
}
\caption{Example 3: Total number of iterations and condition numbers for fixed $ dx=1/10 $, and different $ \Delta t $. Here, $L_1 = L_2 = L_3 = 1$, different $m$ and $m_{lin}$ are taken into account}
\label{tab:3.0dt}
\end{table}

Once more, we report the numerical errors and the estimated orders of convergence associated with the discretization technique here implemented (TPFA). In Table \ref{tab:3EOC}, we present the values obtained for the monolithic L-scheme. The EOC depends only on the discretization technique, not the linearization scheme or solving algorithm. 

\begin{table}[h!]
\begin{center}
	\captionsetup{justification=centering,margin=2cm}
\begin{tabular}{|l|l|l|l|l|l|}
\hline
      &$e_1$ & EOC & $e_2$ & EOC & $e_3$ \\ \hline
 $\Psi$   &   0.0308  &   1.0345        &       0.0150        &  0.9281     &    0.0085       \\ \hline
 $\theta$ &  0.0350     &   1.2985     &    0.0142       &     1.0997     &   0.0066          \\ \hline
 $c$     &   0.0060    &   1.3202    &       0.0024     &     1.3000      &      0.0010      \\ \hline
       \end{tabular}
\caption{Example 3: Numerical error and estimated order of convergence (EOC) of the discretization method.}
\label{tab:3EOC}
\end{center}
\end{table}

\begin{figure}[h!]
	\begin{center} 
		\captionsetup{justification=centering,margin=2cm}
	  \begin{subfigure}{.49\textwidth}
			\includegraphics[width=1\linewidth]{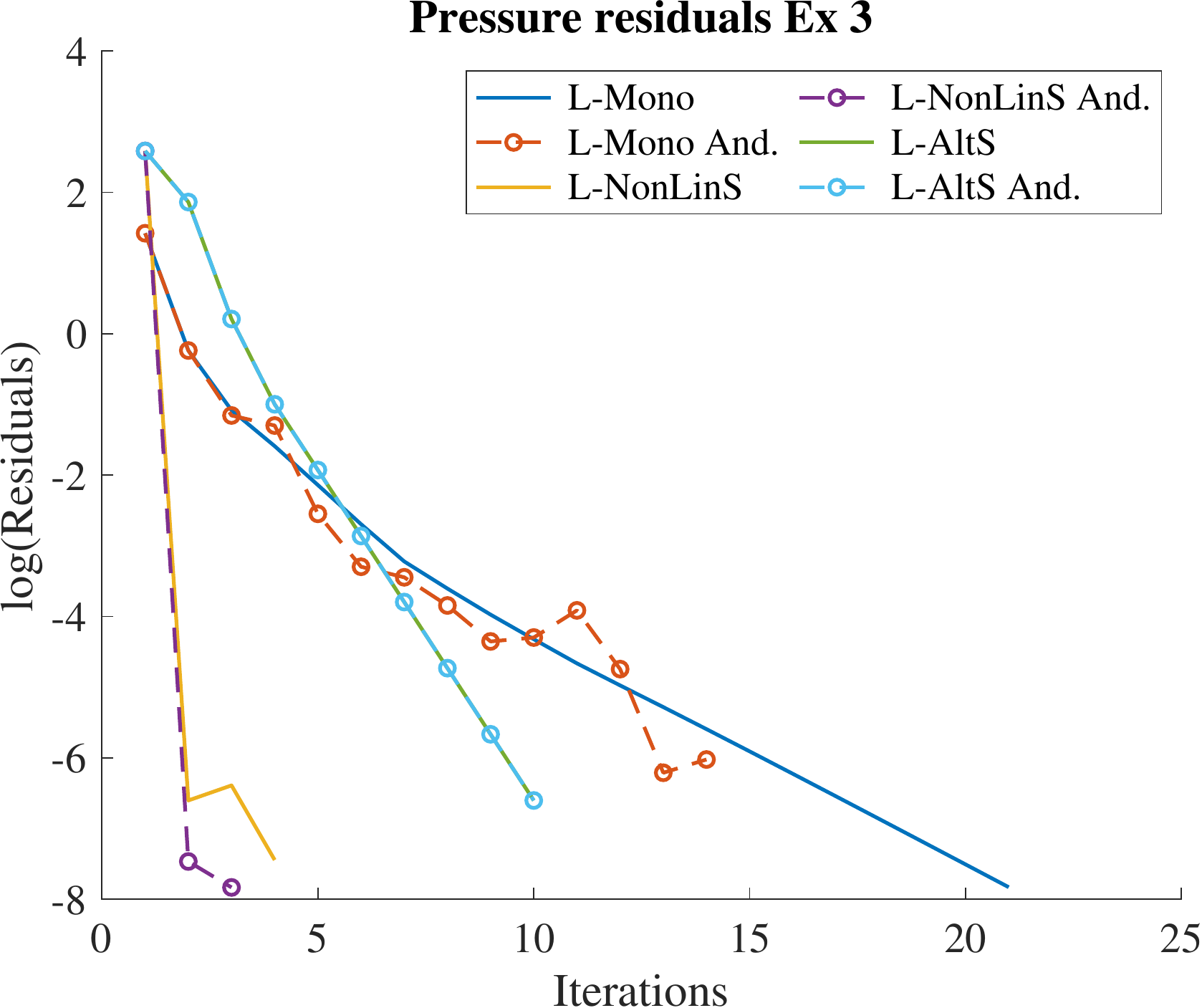}
			\caption{Pressure residuals.}
			\end{subfigure}
  \begin{subfigure}{.49\textwidth}
		  \includegraphics[width=1\linewidth]{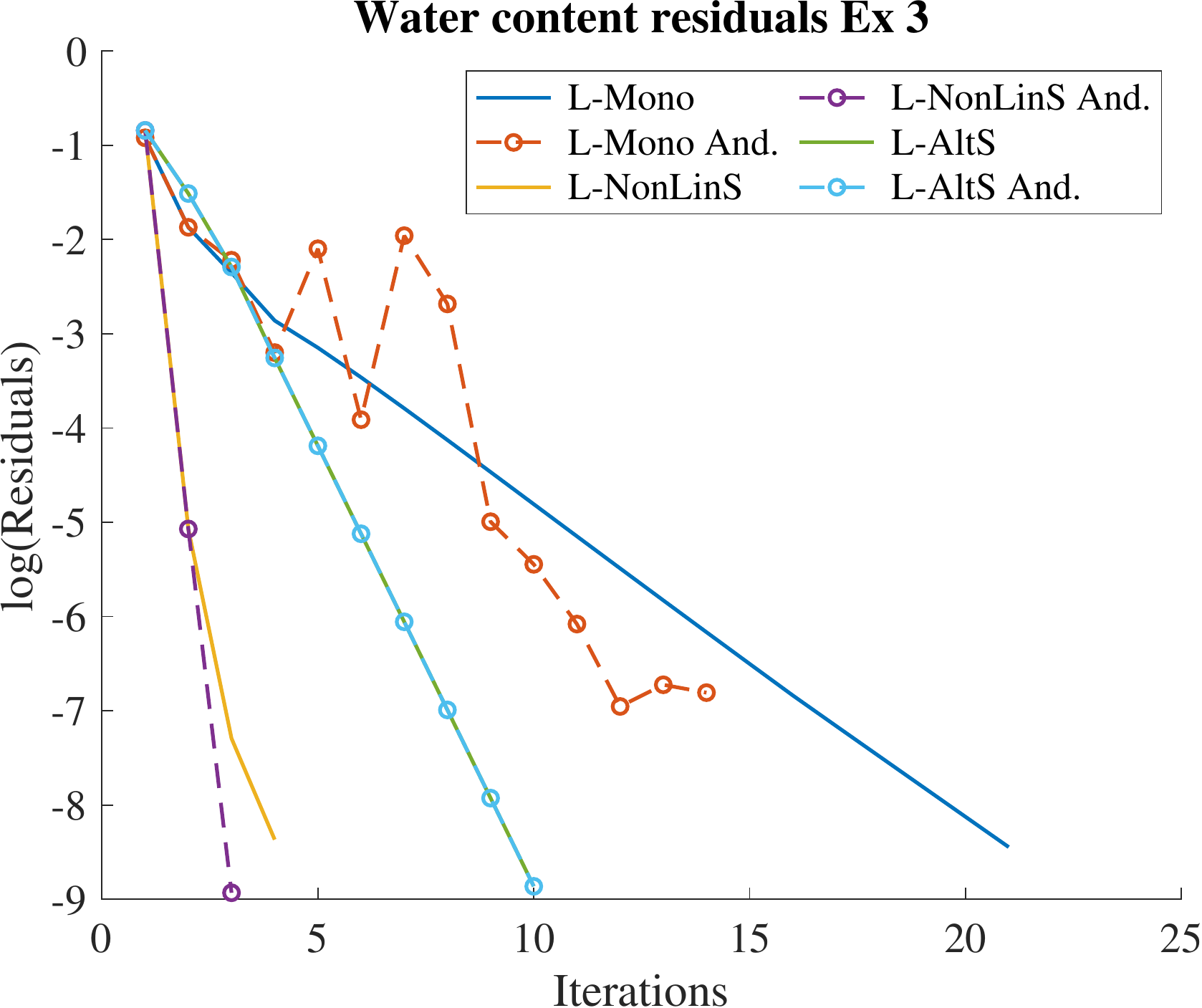}
			\caption{Water content residuals.}
			\end{subfigure}\\
			  \begin{subfigure}{.49\textwidth}
\includegraphics[width=1\linewidth]{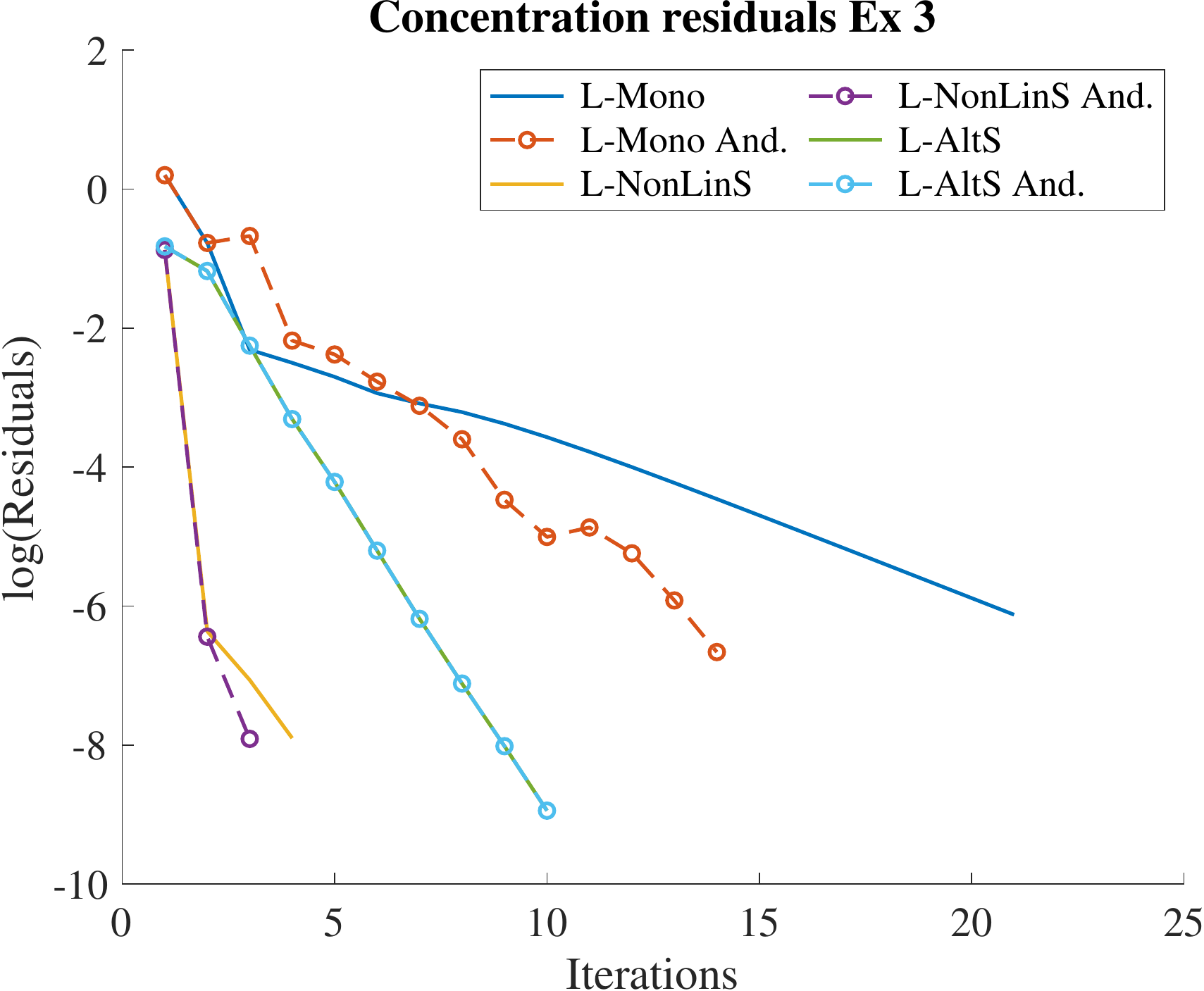}
			\caption{Concentration residuals.}
			\end{subfigure}		
		\caption{Example 3: residuals of each unknown at the final time step, for the different schemes. \\Here, $L_1=L_2=L_3=1$, $dx=1/40$, $\Delta t=T/25$ and $ m \neq m_{lin} $.}
		\label{fig:Ex3Residuals}
	\end{center}
\end{figure}

In Figure \ref{fig:Ex3Residuals}, we report the residuals of pressure, water content and concentration, at the final time step. The differences between the accelerated and non-accelerated schemes seem to be minimal at the final time step but we observe in Tables \ref{tab:3.0dx} and \ref{tab:3.0dt} that the total improvements are actually substantial. The precise orders of convergence for the different solving algorithms are reported in Table \ref{tab:3.1EOC}. The non-accelerated L-schemes have an order of convergence equal to one, while the accelerated ones have slightly larger values. No result was reported for the Newton schemes due to the lack of convergence.

\begin{table}[h!]
\begin{center}
	\captionsetup{justification=centering,margin=2cm}
\begin{tabular}{|c|c|c|c|c|c|c|}
\hline
&\scriptsize{LS-Mono} & \scriptsize{LS-Mono And.} & \scriptsize{LS-NonLinS} & \scriptsize{LS-NonLinS And.} & \scriptsize{LS-AltLinS} & \scriptsize{LS-AltLinS And.}\\ \hline
$\Psi$   &  1.14   &  1.33  & 1.00  &  1.19   &  1.00    &   1.81      \\ \hline
$c$        &  1.00   &  1.32  & 1.00  &  1.51   &   0.99   &   1.52      \\ \hline
$\theta$ &  1.00   &  1.29 & 1.00  &  1.39   &  1.00    &  1.16       \\ \hline
\end{tabular}
\caption{Example 3: Estimated order of convergence for the different linearization schemes.}
\label{tab:3.1EOC}
\end{center}
\end{table}

\newpage
\subsection{Example 4, $\gamma=1$, $\delta=5e-3$, $\tau(\theta)=1+ \theta^2$}

Finally, we study a problem which includes both hysteresis and dynamic capillary effects. We choose $\delta=5e-3$, $\gamma = 1$ and $\tau(\theta)=1+ \theta^2$. As for the previous examples, we report the total numbers of iterations required by each algorithm, the condition numbers associated with the linearized equations, the EOC of the discretization technique and the residual for each unknown.

In Tables \ref{tab:4.0dx} and \ref{tab:4.0dt}, we present the total number of iterations required by each algorithm, and the condition numbers associated with each system.
As in the previous example, the Newton method fails to converge, while the L-scheme based solvers present no difficulties. The L parameters are all set equal to $0.1$. This leads to a faster convergence, when compared to the previous example, where larger values have been required for robustness. This is explained by the fact that, since the dynamic effects are introduced ($ \tau>0 $), the solution is more regular \cite{Cao2015,Mikelic2009}.

We have tested different values of $m$ on the Newton methods, but none ensured the convergence of the schemes. As in the previous test cases, we have investigated smaller time steps, but the Newton solvers has still failed to converge. 

Regarding the results obtained thanks to the AA, we can notice some improvements which are smaller than the ones observed for the previous test cases. Once more, this is due to the optimal choice of the $L$ parameters, ensuring that the L-scheme converges, on average, in 5 iterations per time step. Therefore further improvements are not expected. Note that, compared to the first example (Table \ref{tab:MsEx1}), larger $L$ values are used leading to larger numbers of iterations. This explains why the AA with proper parameters $ m $ have improved the convergence behaviour of the L-scheme there.

\begin{table}[H]
\centering
\captionsetup{justification=centering,margin=2cm}
    \scalebox{.58}{
  \begin{tabular}{ c | c c | c c  c c  | c c c c}
\hline \hline
& Monolithic & &  & NonLinS &  & & & AltLinS &&\\
\hline
& Newton & & & Newton &  & & &Newton  & &\\
\hline
&  & & &  &  cond. \# & & &  & cond. \# &\\
\hline
dx &\# iterations & condition \#  & \# iterations & Flow && Transport  & \# iterations & Flow && Transport\\ 
1/10  & - & - & - & - && -& - & - && -\\ 
1/20  & - & - & - & - && -& - & - && - \\ 
1/40  & - & - & - & - && -& - & - && - \\  
 \hline \hline
 & L-scheme & & & L-scheme &  & & & L-scheme & & \\
 \hline
 &  & & &  &  cond. \# & & &  & cond. \# &\\
\hline
dx &\# iterations & condition \#  & \# iterations & Flow && Transport  & \# iterations & Flow && Transport\\ 
1/10  & 152   & 290.13  & 128 - 122  & 208.09  && 162.46 & 251   &   206.24 && 174.77  \\ 
1/20  &  160  &  768.54  &  137 - 121 &  558.01 && 621.77 & 259   &  486.23  && 598.01 \\ 
1/40  &  165  &  3.04e+03  & 141 - 120  & 2.10e+03  && 2.35e+03 &  328  & 2.14e+03 && 2.4415e+03\\ 
 \hline \hline
 & L-scheme & (AA m = 1)  & & L-scheme (AA m = 2, $m_{lin}$ = 3) &  & & & L-scheme (AA m = 2) & & \\
 \hline
 &  & & &  &  cond. \# & & &  & cond. \# &\\
\hline
dx &\# iterations & condition \#  & \# iterations & Flow && Transport  & \# iterations & Flow && Transport\\ 
1/10  & 139  & 288.40  & 112 - 85  & 212.29  && 165.24 &152  &198.15   && 166.47 \\ 
1/20  &  144 & 752.41  & 117 - 89  & 550.05  && 630.43 & 162  &522.01   &&  636.00\\ 
1/40  &  149 & 3.04e+03& 127 - 88  & 2.19e+03  && 2.39e+03 &  166  & 2.06e+03  && 2.40e+03 \\ 
 \hline \hline
\end{tabular}
}
\caption{Example 4: Total number of iterations and condition numbers for fixed $\Delta t=T/25$, and different $ dx $. Here, $L_1 = L_2 = L_3 = 0.1$ and $m\neq m_{lin}$.}
\label{tab:4.0dx}
\end{table}

\begin{table}[H]
\centering
\captionsetup{justification=centering,margin=2cm}
    \scalebox{.58}{
  \begin{tabular}{ c | c c | c c  c c  | c c c c}
\hline \hline
& Monolithic & &  & NonLinS &  & & & AltLinS &&\\
\hline
& Newton & & & Newton &  & & &Newton  & &\\
\hline
&  & & &  &  cond. \# & & &  & cond. \# &\\
\hline
$\Delta t$ &\# iterations & condition \#  & \# iterations & Flow && Transport  & \# iterations & Flow && Transport\\ 
T/25   & - & - & - & -  && - & - & - && -\\ 
T/50   & - & - & - & -  && - & - & - && -\\ 
T/100 & - & - & - & -  && - & - & - && -\\
 \hline \hline
 & L-scheme & & & L-scheme &  & & & L-scheme & & \\
 \hline
&  & & &  &  cond. \# & & &  & cond. \# &\\
\hline
$\Delta t$ &\# iterations & condition \#  & \# iterations & Flow && Transport  & \# iterations & Flow && Transport\\ 
T/25    & 152   & 290.13  & 128 - 122  & 208.09  && 162.46 & 251   &   206.24 && 174.77  \\ 
T/50    &  248  & 310.93  &  201 - 225 & 264.09 &&  89.64  &  424   &  263.91  && 97.09  \\ 
T/100  &  508  & 415.80  &  403- 405   & 403.31  && 47.53  &  768   &  403.81  && 52.72  \\ 
 \hline \hline
 & L-scheme & (AA m = 1)  & & L-scheme (AA m = 2, $m_{lin}$ = 3) &  & & & L-scheme (AA m = 2) & & \\
 \hline
 &  & & &  &  cond. \# & & &  & cond. \# &\\
\hline
$\Delta t$ &\# iterations & condition \#  & \# iterations & Flow && Transport  & \# iterations & Flow && Transport\\ 
T/25    & 139  & 288.40 & 112 - 85  & 212.29  && 165.24 & 152  &198.15   && 166.47 \\ 
T/50    & 233  & 312.95 &  195 - 167   &   267.27&&  90.55  &    250   &  260.41 &&  91.42 \\
T/100  & 448  &  416.05&  358 - 308   & 404.34  &&  47.65  &  506    &   403.30   &&   48.38  \\
 \hline \hline
\end{tabular}
}
\caption{Example 4: Total number of iterations and condition numbers for fixed $dx=1/10$, and different $ \Delta t $. Here  $L_1 = L_2 = L_3 = 0.1$ and $m\neq m_{lin}$ }
\label{tab:4.0dt}
\end{table}

The numerical errors and the estimated orders of convergence of the discretization technique (TPFA), presented in Table \ref{tab:4EOC}, are consistent with the ones from the previous test cases.

\begin{table}[h!]
\begin{center}
	\captionsetup{justification=centering,margin=2cm}
\begin{tabular}{|l|l|l|l|l|l|}
\hline
&$e_1$ & EOC & $e_2$ & EOC & $e_3$ \\ \hline
 $\Psi$   &   0.0759  & 0.9558  &    0.0391  & 0.8837 &   0.0212   \\ \hline
 $\theta$ &0.0138   &0.9115   &   0.0073   &          0.9463 &   0.0038   \\ \hline
 $c$     &0.0101   &1.2531   &   0.0042 &   1.2655 &   0.0018  \\ \hline
       \end{tabular}
\caption{Example 4: Numerical error and estimated order of convergence (EOC) of the discretization method.}
\label{tab:4EOC}
\end{center}
\end{table}

\begin{figure}[h!]
	\begin{center} 
		\captionsetup{justification=centering,margin=2cm}
	  \begin{subfigure}{.49\textwidth}
			\includegraphics[width=1\linewidth]{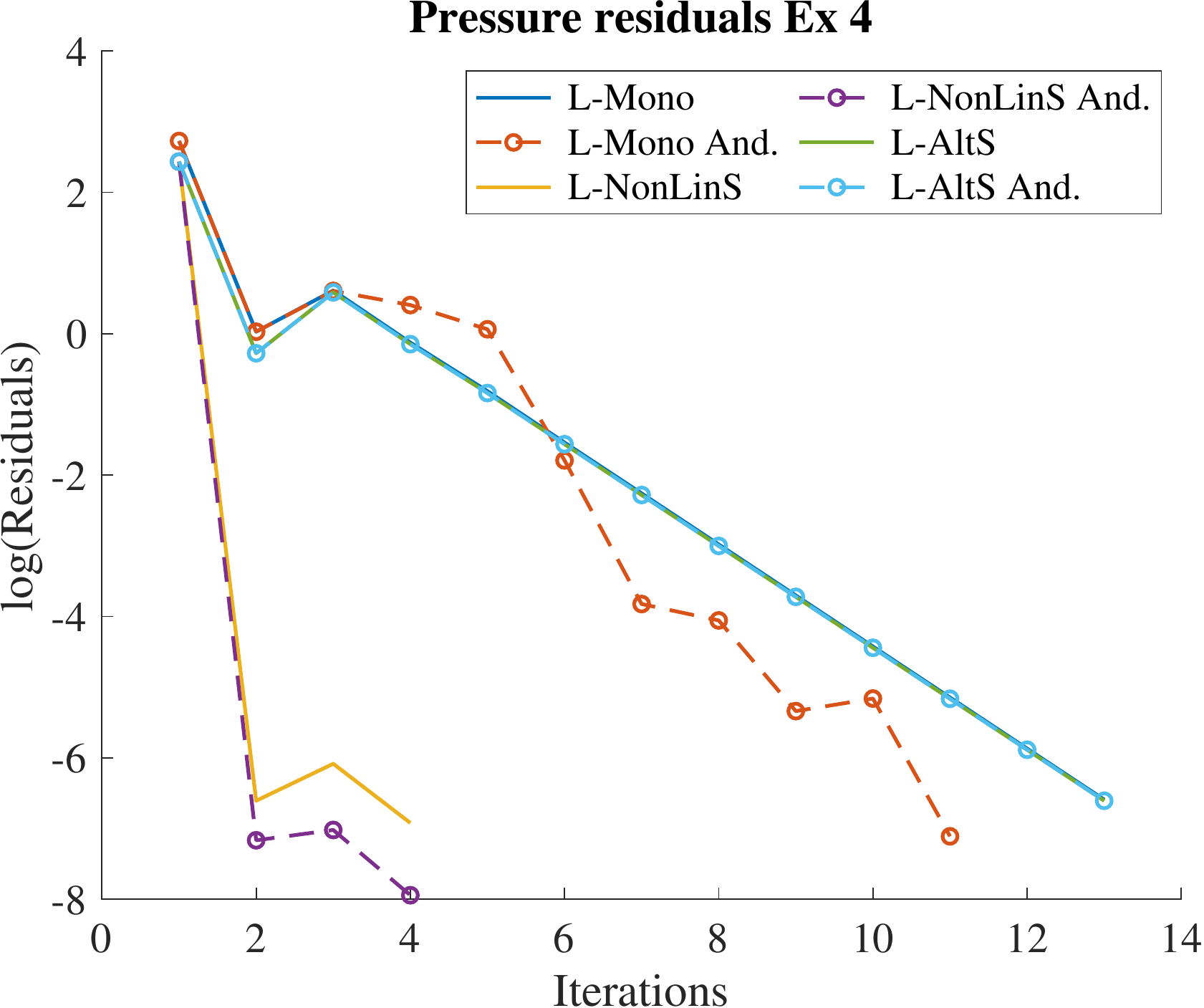}
			\caption{Pressure residuals.}
			\end{subfigure}
  \begin{subfigure}{.49\textwidth}
		  \includegraphics[width=1\linewidth]{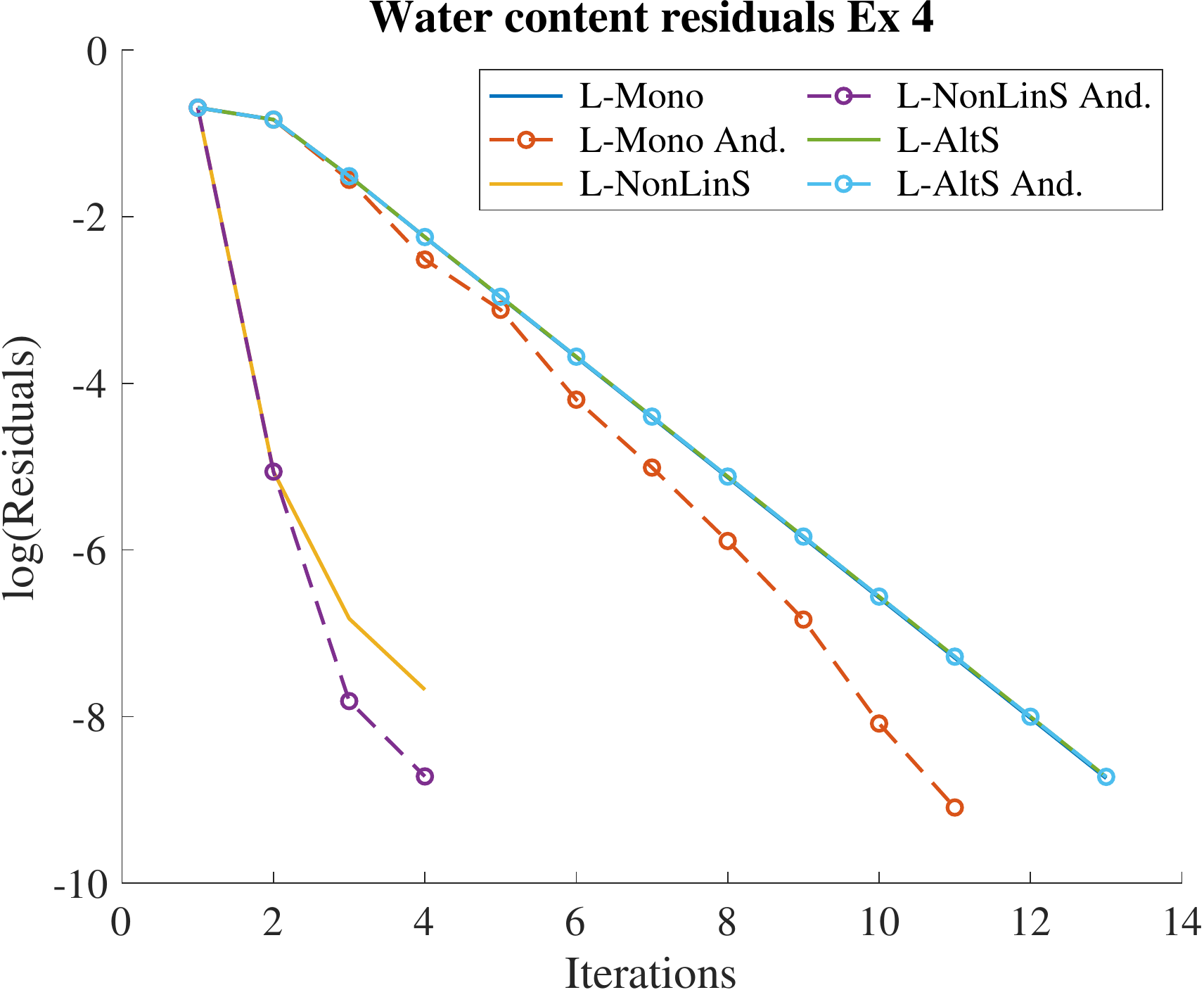}
			\caption{Water content residuals.}
			\end{subfigure}\\
			  \begin{subfigure}{.49\textwidth}
\includegraphics[width=1\linewidth]{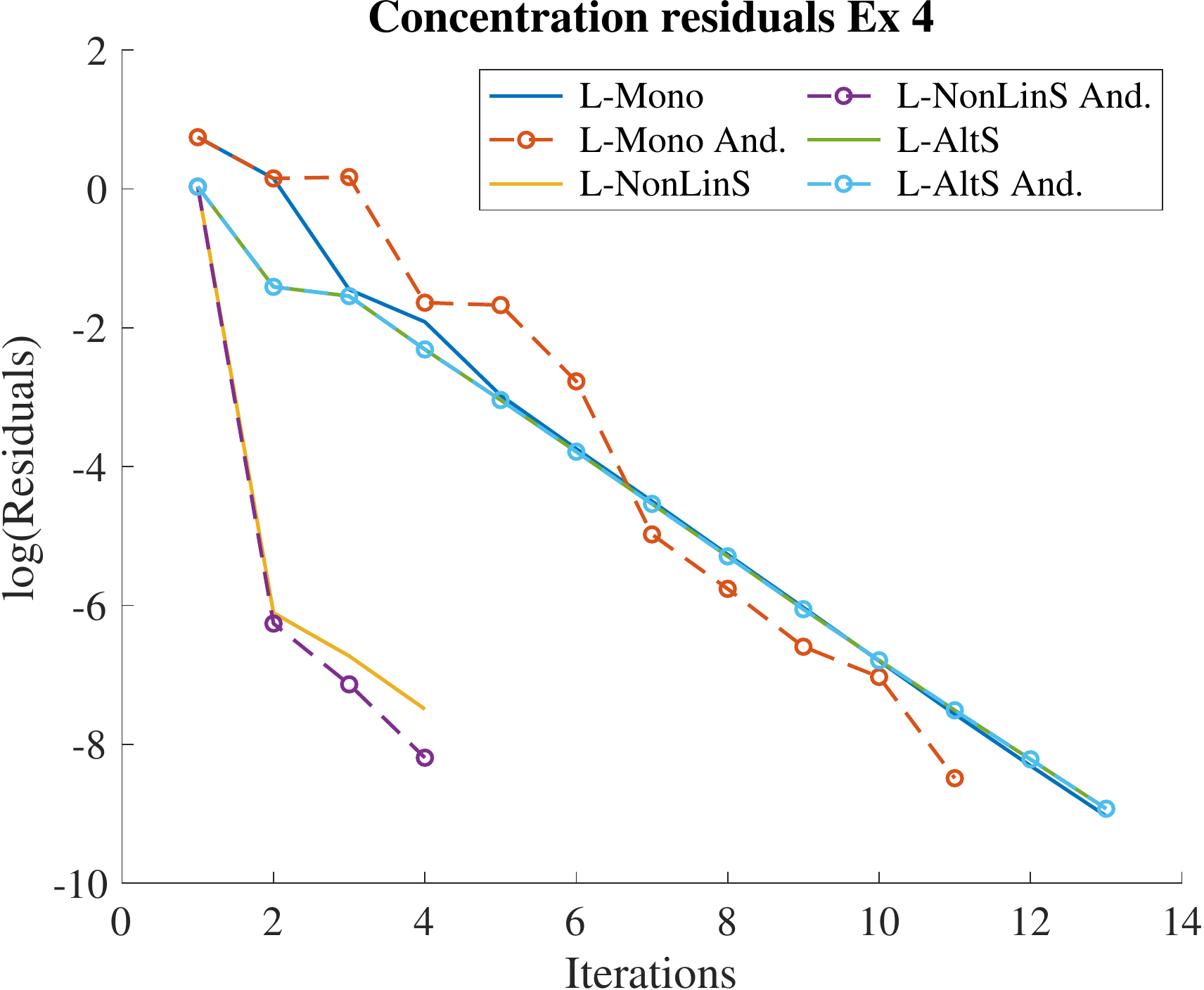}
			\caption{Concentration residuals.}
			\end{subfigure}		
		\caption{Example 4: Residuals of each unknown at the final time step, for the different schemes. \\Here, $L_1=L_2=L_3=0.1$, $dx=1/40$, $\Delta t=T/25$ and different $ m \neq m_{lin} $.}
		\label{fig:Ex4Residuals}
	\end{center}
\end{figure}

Finally, regarding the order of convergence of the different solving algorithms, Figure \ref{fig:Ex4Residuals} presents the residuals for each unknown, and Table \ref{tab:4.1EOC} the precise orders computed averaging over iterations at the final time step. 

\begin{table}[h!]
\begin{center}
	\captionsetup{justification=centering,margin=2cm}
\begin{tabular}{|c|c|c|c|c|c|c|}
\hline
&\scriptsize{LS-Mono} & \scriptsize{LS-Mono And.} & \scriptsize{LS-NonLinS} & \scriptsize{LS-NonLinS And.} & \scriptsize{LS-AltLinS} & \scriptsize{LS-AltLinS And.}\\ \hline
$\Psi$   & 1.00  & 1.29  & 0.99    &  1.27 & 0.99 & 1.30    \\ \hline
$c$        & 1.29 &  1.29 &  1.02   &  1.15 &  1.07&  1.11   \\ \hline
$\theta$ & 1.00 &  1.16& 1.00    &  1.35 & 1.00 &  1.20   \\ \hline
\end{tabular}
\caption{Example 4: Estimated order of convergence for the different linearization schemes.}
\label{tab:4.1EOC}
\end{center}
\end{table}

\newpage
\subsection{Physical example}

As final numerical study, we investigate a test case that involves realistic parameters, but without having a manufactured solution. The flow will be driven by the boundary conditions. The domain $\Omega$ is the unit square and the final time is $T= 4$. The capillary pressure and conductivity expressions are given by the van Genuchten formulation \cite{vanG}, $K(\theta) = \theta_e^l \Big(1-\big(1-\theta_e^{1/M}\big)^M\Big)^2$ and $p_{cap}(\theta,c)=\big(1-b \ln(c/a + 1)\big)^{-1}(-\theta^{-1/M})^{1-M}$, \\where $\theta_e = (\theta-\theta_r)/(\theta_s-\theta_r)$ is the effective water content, $\theta_s=0.9$, $\theta_r=0.005$, $M=2$, $l= 0.31$, $a=0.04$ and $b=0.47$. Furthermore, we take $\tau(\theta)=1 + \theta^2$, and the hysteresis effects are included by setting $\gamma = 1$ and $\delta = 5e-3$, as in Example 4.

Dirichlet boundary conditions are imposed at the left side of the unit square
\begin{equation*}
\begin{split}
 \Psi|_{x=0} &= 1 + \begin{cases}  0.5t \qquad &\text{if}\  t<1, \\
  0.5 \qquad &\text{if}\  1\leq t<2, \\
0.5(3-t)\qquad &\text{if}\  2\leq t<3, \\
 -0.4\qquad &\text{if}\  3\leq t\leq4, \\
\end{cases}\\
c|_{x=0} &= 2,
\end{split}
\end{equation*}
whereas, at the remaining sides, we consider homogeneous Neumann boundary conditions. The discontinuity in time $t=3$ makes solving the problem numerically even more complex. The initial conditions are
\begin{equation*}
\theta^0 = x, \qquad \text{and} \qquad c^0 = 1.
\end{equation*}

All L parameters are set to $0.5$. We have tested different values, but, $L=0.5$ seems to give the best results in terms of numbers of iterations. Furthermore, the results may be improved even further by choosing different values for each parameter, $L_1,L_2$ and $L_3$, but for ease of presentation this has been omitted here.

In Tables \ref{tab:5.0dx} and \ref{tab:5.0dt}, we report the total numbers of iterations and condition numbers associated to each algorithm. We observe that, due to the higher nonlinearities of the conductivity $K$ and capillary pressure $p_{cap}$ involved, the results are different compared to the ones presented for the previous examples. 

Again, the Newton solvers have failed to converge and the systems associated with the linearized equations are badly conditioned. Considering smaller time steps did not resolve this.

The L-schemes on the other hand converge, but require high numbers of iterations. We observe a significant improvement thanks to the AA. The performance of the monolithic solver is for example drastically improved, the iterations required are reduced by circa 50\%.  In case of finer meshes, one needs to use a larger L, precisely $L=1$. The AA can introduce some instabilities, and thus a larger L may be required. Clearly, this leads to an increase in the number of iterations. Such results are still better than the one obtained for smaller L without acceleration. Similar observations can also be made for the splitting solvers. Even though larger L parameters may be required, the accelerated schemes perform better then the non-accelerated ones. The nonlinear splitting seems to be less stable than the alternate linear one. For a coarse mesh, we could use a large $m$, $m=5$, resulting in an extreme reduction in the numbers of iterations. For finer meshes, $m$ had to be set equal to 1, otherwise the schemes did not converge, and a larger L parameter was required.
The alternate linearized splitting seems to be more stable, as the L-scheme linearization acts as stabilization term for both the nonlinearities and the decoupling. Unfortunately, for the first time, the results are slower than for the nonlinear splitting. The main differences are observable in the second table, Table \ref{tab:5.0dt}. The mesh size is fixed, $dx = 1/10$, and the nonlinear solvers converge even for $m=5$. The large AA parameter $m$ ensures a strong reduction in the numbers of iterations. The alternate linearized splitting converges only for $m=1$, anyhow the improvements are remarkable.

\begin{table}[H]
\centering
\captionsetup{justification=centering,margin=2cm}
    \scalebox{.5}{
  \begin{tabular}{ c | c c | c c  c c  | c c c c}
\hline \hline
& Monolithic & &  & NonLinS &  & & & AltLinS &&\\
\hline
& Newton & & & Newton &  & & &Newton  & &\\
\hline
&  & & &  &  cond. \# & & &  & cond. \# &\\
\hline
dx &\# iterations & condition \#  & \# iterations & Flow && Transport  & \# iterations & Flow && Transport\\ 
1/10  & - & - & - & - && -& - & - && -\\ 
1/20  & - & - & - & - && -& - & - && - \\ 
1/40  & - & - & - & - && -& - & - && - \\  
 \hline \hline
 & L-scheme & & & L-scheme &  & & & L-scheme & & \\
 \hline
 &  & & &  &  cond. \# & & &  & cond. \# &\\
\hline
dx &\# iterations & condition \#  & \# iterations & Flow && Transport  & \# iterations & Flow && Transport\\ 
1/10  &  899   & 702.16     & 1522 - 235   &  790.37  &&  155.18   &    1490 &  623.78  &&  155.10  \\
1/20  &  930   & 3.58e+03 &  1515 - 240  & 3.32e+03  && 615.72 &    1428 &  2.54e+03  &&  596.60  \\
1/40  &  941   & 1.52e+04 &  1680 - 243  & 1.44e+04   &&  2.46e+03   &  1548   &  1.24e+04  &&  2.45e+03  \\
 \hline \hline
 & L-scheme & (AA m = 1)  & & L-scheme  (AA m = 5, $m_{lin}$ = 1) &  & & & L-scheme (AA m = 1) & & \\
 \hline
 &  & & &  &  cond. \# & & &  & cond. \# &\\
\hline
dx &\# iterations & condition \#  & \# iterations & Flow && Transport  & \# iterations & Flow && Transport\\ 
1/10  &  480   & 1.37e+03  &  210 - 96  & 811.29 && 155.16  & 795    &  774.00  &&  155.18  \\
1/20  &  541 ($L=1$)   & 5.27e+03  & 758 - 170 ($m=1, L=1$)&   2.76e+03 &&  395.74   &  532   &  2.72e+03  &&  468.04  \\
1/40  &  603 ($L=1$)   &  2.29e+04 & 1214 - 261 ($m=1, L=2$)  &   8.84e+03 && 911.31 & 1798 ($L_2=2$)  & 1.89e+04   &&2.46e+03  \\
\hline\hline
\end{tabular}
}
\caption{Example 5: Total number of iterations and condition numbers for fixed $\Delta t=T/25$, and different $ dx $. Here, $L_1 = L_2 = L_3 = 0.5$ and $m \neq m_{lin}$.}
\label{tab:5.0dx}
\end{table}

\begin{table}[H]
\centering
\captionsetup{justification=centering,margin=2cm}
    \scalebox{.55}{
  \begin{tabular}{ c | c c | c c  c c  | c c c c}
  \hline\hline
& Monolithic & &  & NonLinS &  & & & AltLinS &&\\
\hline
& Newton & & & Newton &  & & &Newton  & &\\
\hline
&  & & &  &  cond. \# & & &  & cond. \# &\\
\hline
dx &\# iterations & condition \#  & \# iterations & Flow && Transport  & \# iterations & Flow && Transport\\ 
T/25   & -  &  - & - & - && -& - & - && -\\
T/50   & -  &  - & - & - && -& - & - && -\\
T/100 & -  &  - & - & - && -& - & - && -\\
 \hline \hline
 & L-scheme & & & L-scheme &  & & & L-scheme & & \\
 \hline
&  & & &  &  cond. \# & & &  & cond. \# &\\
\hline
$\Delta t$ &\# iterations & condition \#  & \# iterations & Flow && Transport  & \# iterations & Flow && Transport\\ 
T/25   & 899   & 702.16  & 1522 - 235   &  790.37  &&  155.18 &    1490 &  623.78  &&  155.10  \\
T/50   & 2892 & 908.39  & 5218 - 450  &  747.79  && 84.75     &    4942 &  533.58  &&  84.83  \\
T/100 & 9809 & 707.12  & 18261 - 849   & 692.77   &&  46.72   &    17406 &  485.71  && 46.84   \\
 \hline \hline
& L-scheme & (AA m =1)  & & L-scheme (AA m = 5, $m_{lin}$ = 1) &  & & & L-scheme (AA m = 1) & & \\
 \hline
 &  & & &  &  cond. \# & & &  & cond. \# &\\
\hline
$\Delta t$ &\# iterations & condition \#  & \# iterations & Flow && Transport  & \# iterations & Flow && Transport\\ 
T/25   & 480 & 1.37e+03  & 210 - 96  & 811.29 && 155.16 &  795    &  774.00  &&  155.18  \\
T/50   & 897 & 1.65e+03  & 517 - 214 & 618.31  &&  53.28 &  1784   &  775.76  &&   84.77 \\
T/100 & 2106 & 1.04e+03& 1354 - 424 & 576.54 && 28.62    &  4928   &  1.05e+03  &&  46.73  \\
\hline\hline
\end{tabular}
}
\caption{Example 5:  Total number of iterations and condition numbers for fixed $ dx=1/20 $, and different $ \Delta t $. Here, $L_1 = L_2 = L_3 = 0.5$ and $m \neq m_{lin}$.}
\label{tab:5.0dt}
\end{table}

It is also interesting to notice that, for decreasing time steps, the number of L iterations per time step is increasing. This is coherent with the theory. With the AA, this is mitigated; the average number of iterations remains more stable.

In Table \ref{tab:5_mValues}, we report the different $m$ and $m_{lin}$ values investigated for the AA. We observe that for the L-scheme Mono solver, the optimal choice, in terms of numbers of iterations is $m=1$. For the nonlinear splitting solver, only one value of $m_{lin}$ has been taken in consideration for the coupling loop, precisely $m_{lin}=1$. This is justified by the fact that the majority of the iterations have taken place in the inside loops, the nonlinearities of the equations are playing a larger role than the coupling aspect.

\begin{table}[H]
\begin{center}
	\captionsetup{justification=centering,margin=2cm}
\begin{tabular}{|l|c|c|c|c|}
\hline
                 				&\# it. m=0 & \# it. m=1 & \# it. m=2 & \# it. m=5 \\ \hline
Mono L-scheme    & 899 &   480 &    775 &   -- \\ \hline
NonLinS L-scheme & 1522 - 325 ($m_{lin}$=0) &  322 -  118 ($m_{lin}$=1)   &  232 - 96  ($m_{lin}$=1) &   210 - 96 ($m_{lin}$=1) \\ \hline
AltLinS  L-scheme  &	   1490   & 795   &   --  &  --  \\ \hline
Newton-Mono        &	  -- &   --  &  -- &   -- \\ \hline
\end{tabular}\caption{Example 5: Numbers of iterations associated to different $m$ and $m_{lin}$ values, $dx=1/10$ $dt=T/25$. }
\label{tab:5_mValues}
\end{center}
\end{table}

\begin{figure}[h!]
	\begin{center} 
		\captionsetup{justification=centering,margin=2cm}
	  \begin{subfigure}{.29\textwidth}
			\includegraphics[width=1\linewidth]{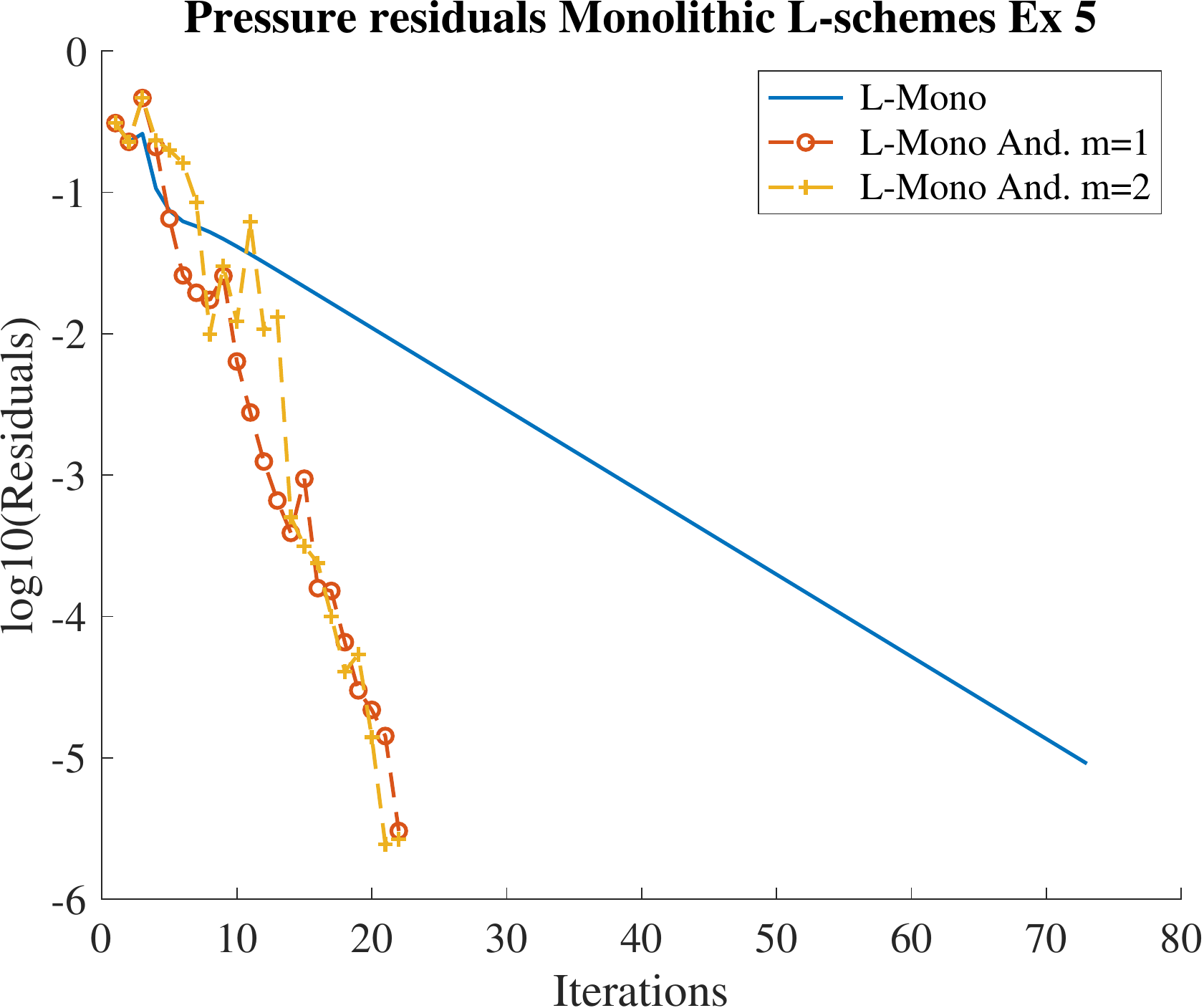}
			\caption{Pressure residuals.}
			\end{subfigure}
  \begin{subfigure}{.29\textwidth}
		  \includegraphics[width=1\linewidth]{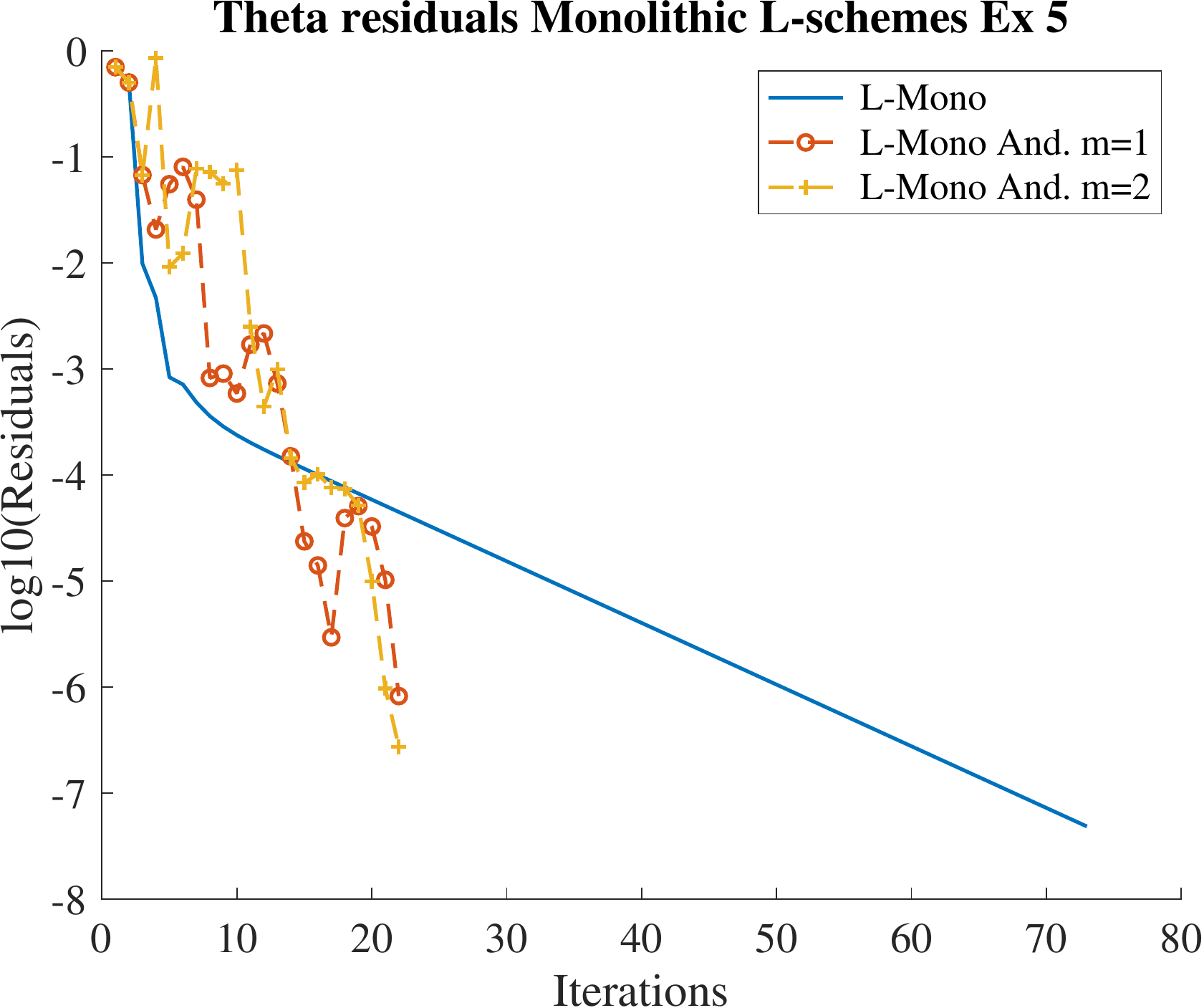}
			\caption{Water content residuals.}
			\end{subfigure}
			  \begin{subfigure}{.29\textwidth}
\includegraphics[width=1\linewidth]{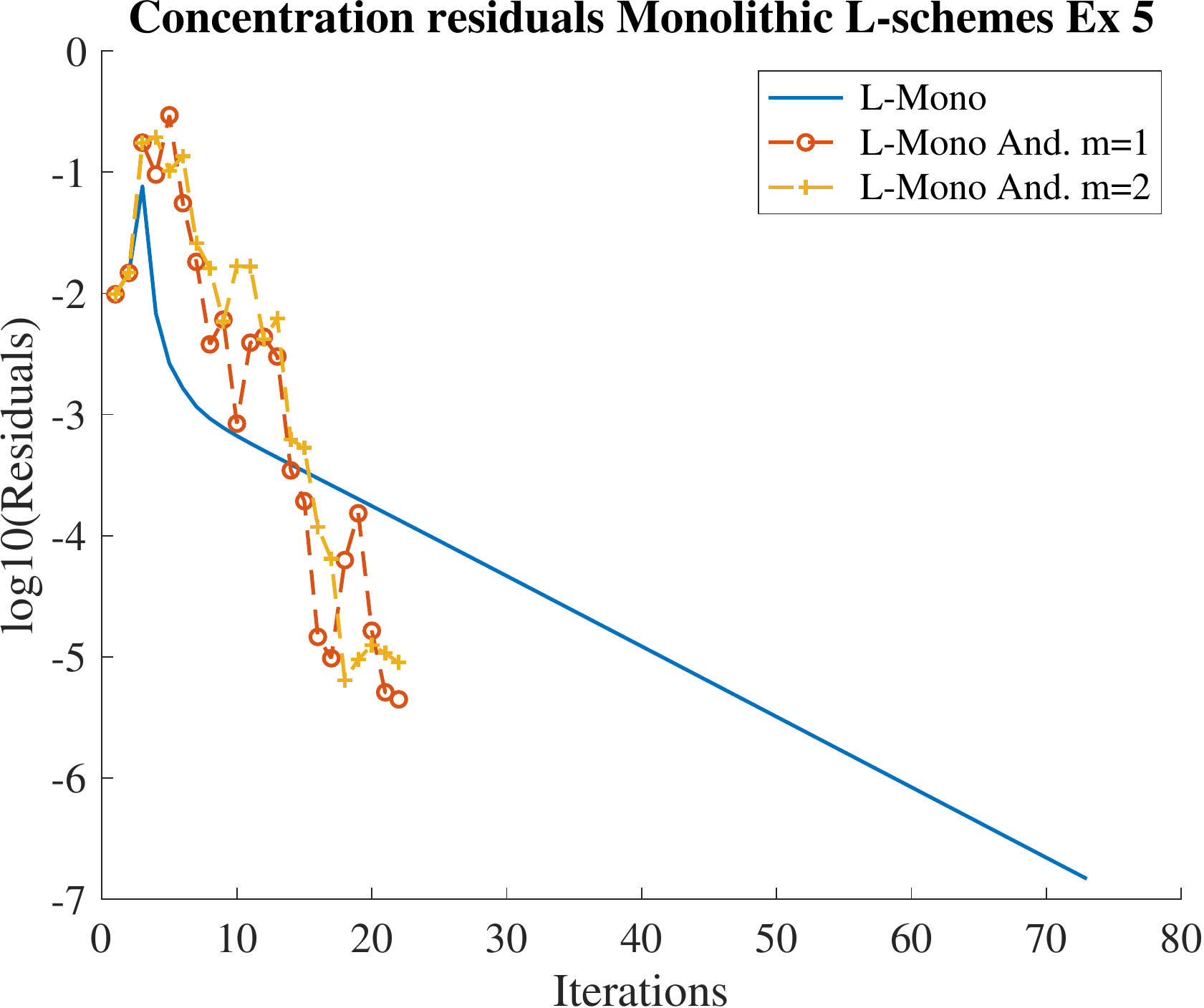}
			\caption{Concentration residuals.}
			\end{subfigure}		
		\caption{Example 5: Residuals of each unknown at the final time step, monolithic L-scheme. Here, $L_1=L_2=L_3=0.1$, different $m$ are tested, $dx=1/10$, and $\Delta t=T/25$.}
		\label{fig:Ex5ResidualsMono}
	\end{center}
\end{figure}

\begin{figure}[h!]
	\begin{center} 
		\captionsetup{justification=centering,margin=2cm}
	  \begin{subfigure}{.29\textwidth}
			\includegraphics[width=1\linewidth]{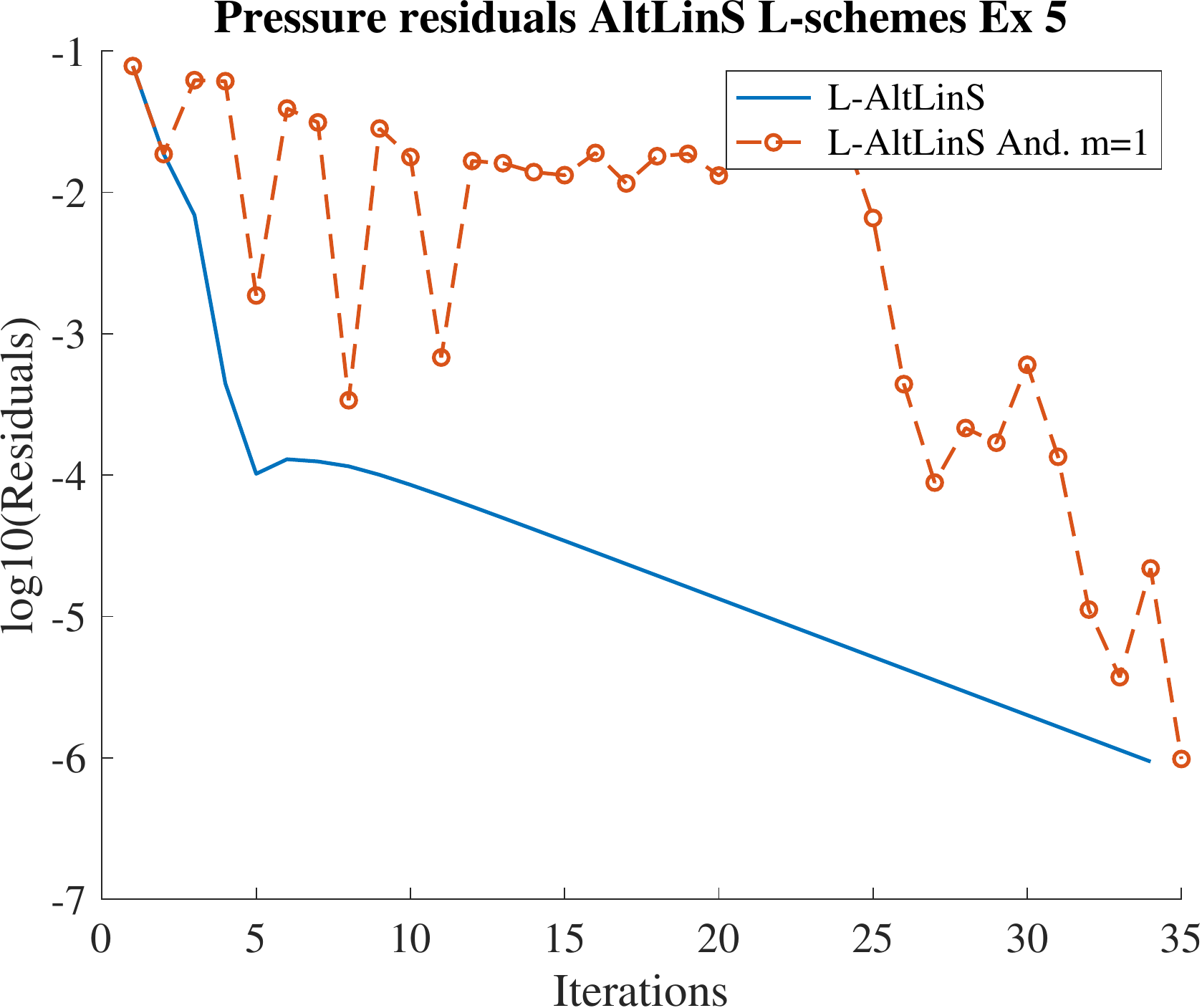}
			\caption{Pressure residuals.}
			\end{subfigure}
  \begin{subfigure}{.29\textwidth}
		  \includegraphics[width=1\linewidth]{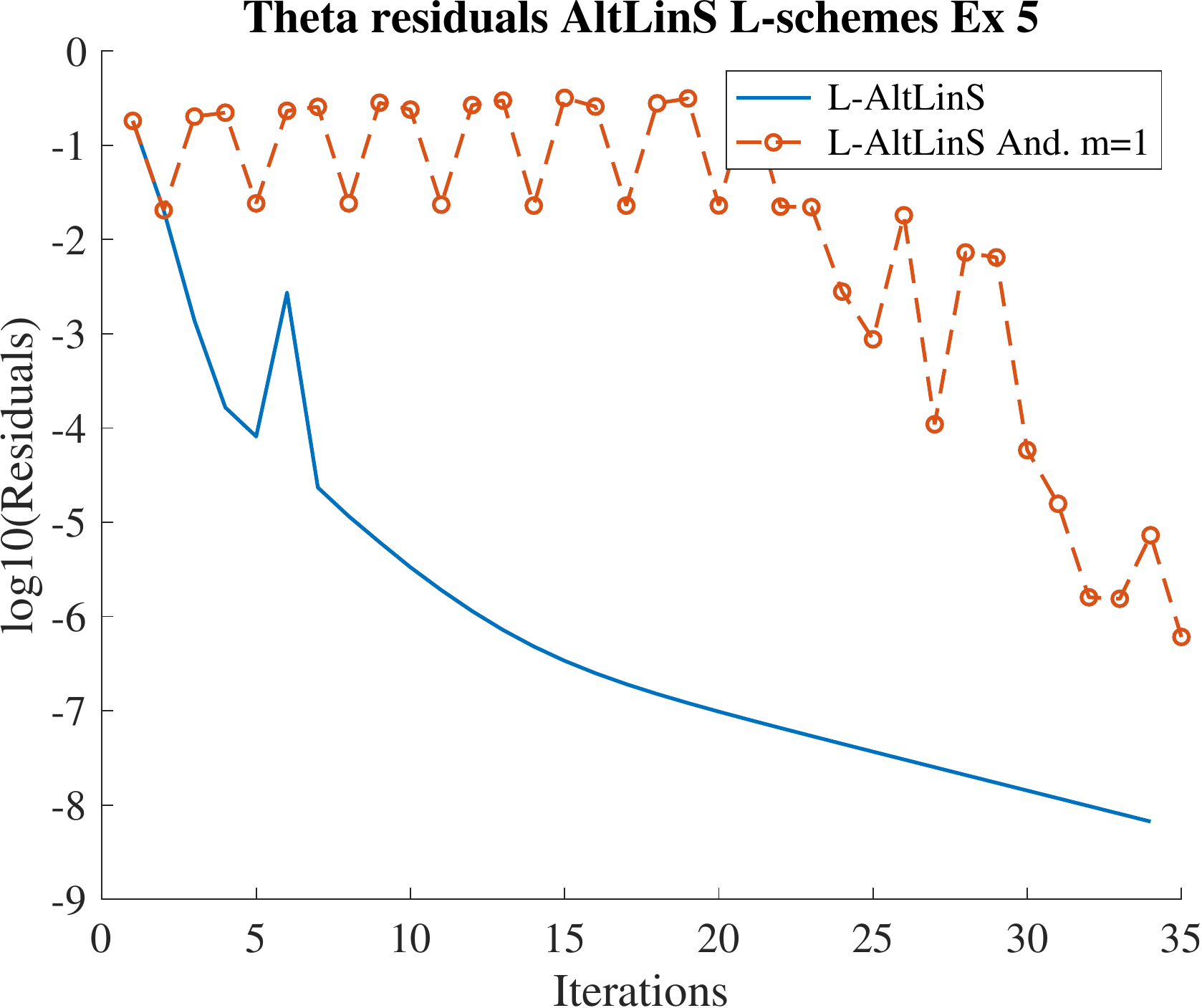}
			\caption{Water content residuals.}
			\end{subfigure}
			  \begin{subfigure}{.29\textwidth}
\includegraphics[width=1\linewidth]{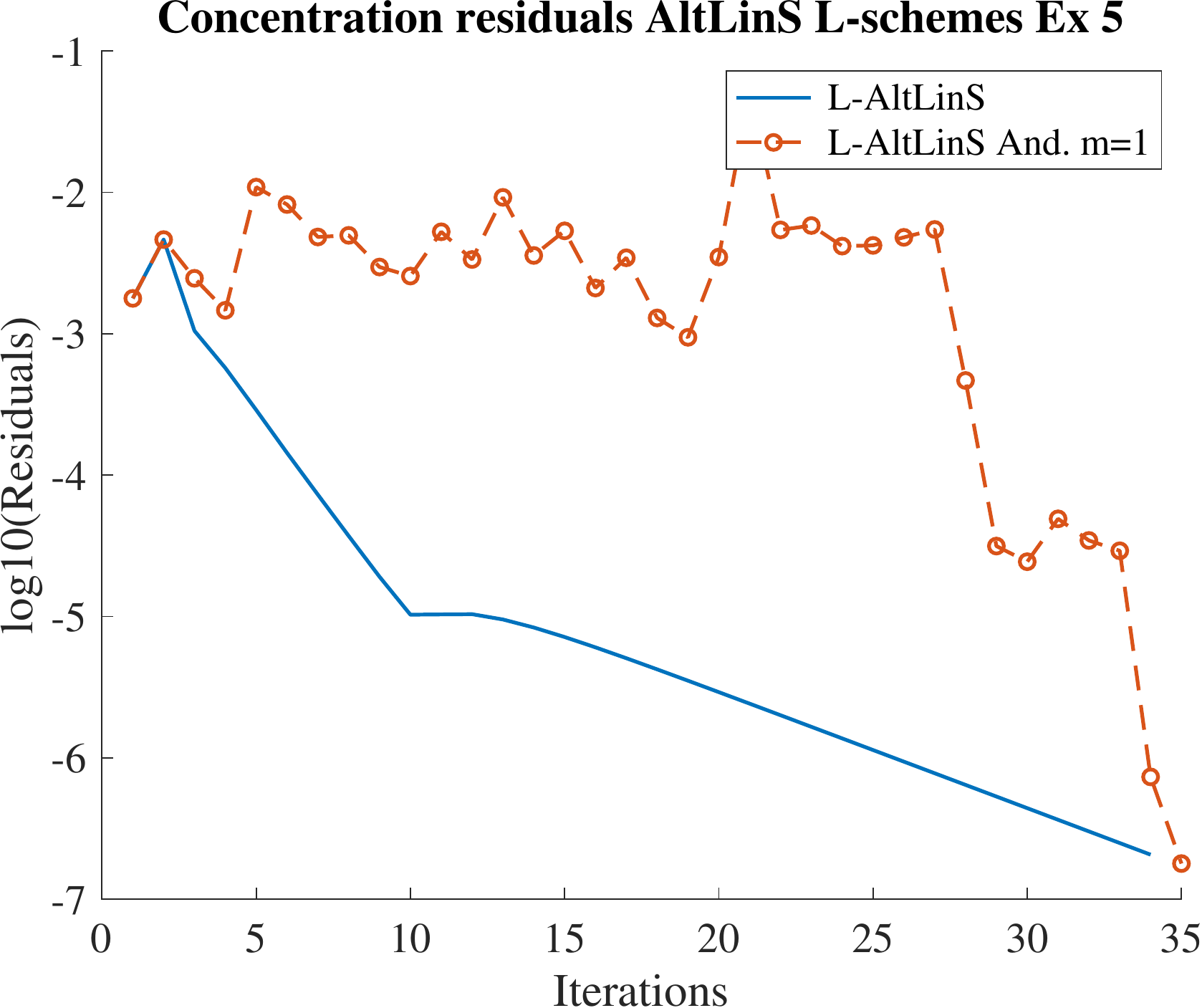}
			\caption{Concentration residuals.}
			\end{subfigure}		
		\caption{Example 5: Residuals of each unknown at the final time step, AltLinS L-scheme. Here, $L_1=L_2=L_3=0.1$, different $m$ are tested, $m_{lin} = 1$, $dx=1/40$, and $\Delta t=T/25$.}
		\label{fig:Ex5ResidualsAltLinS}
	\end{center}
\end{figure}

\begin{figure}[h!]
	\begin{center} 
		\captionsetup{justification=centering,margin=2cm}
	  \begin{subfigure}{.29\textwidth}
			\includegraphics[width=1\linewidth]{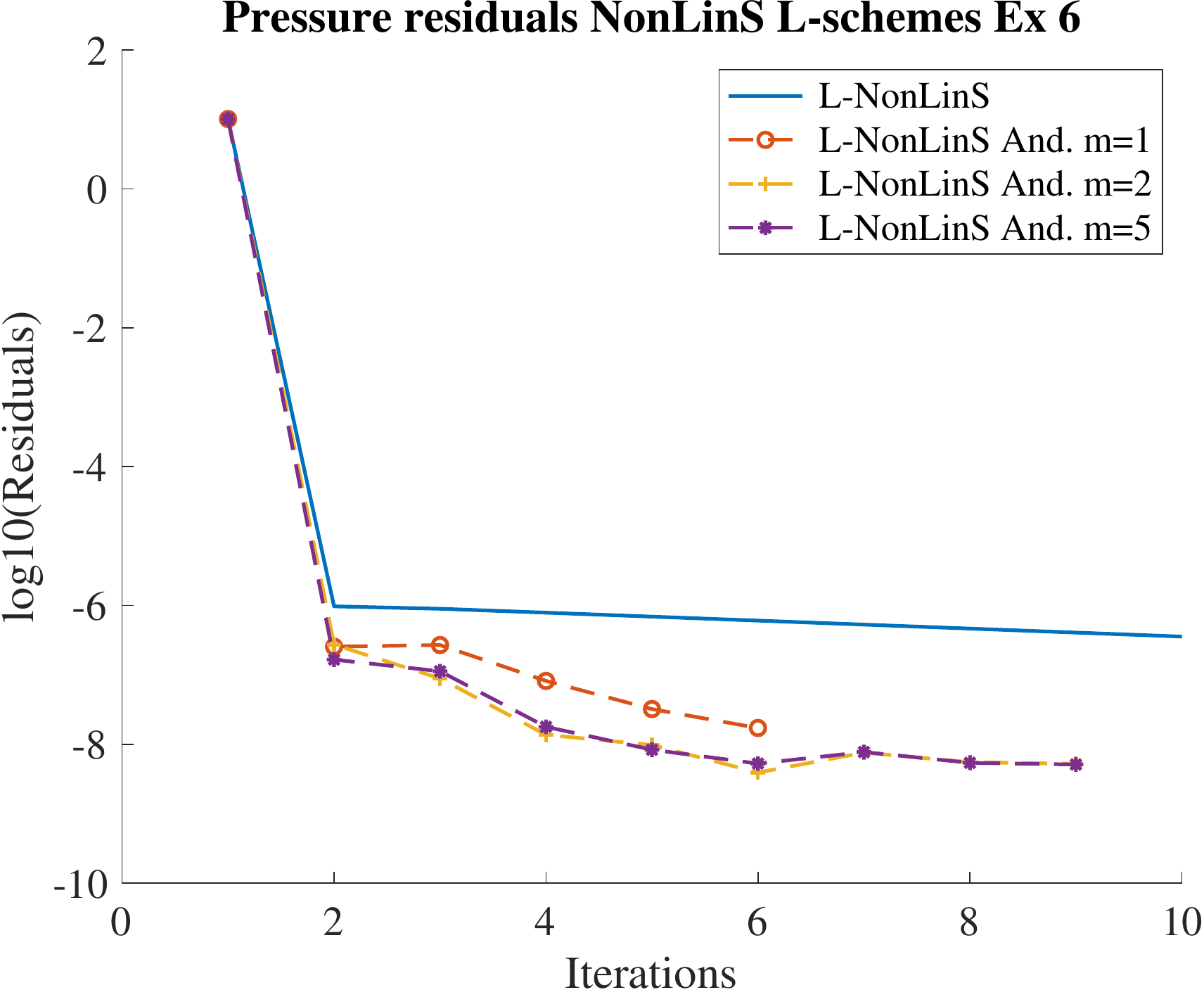}
			\caption{Pressure residuals.}
			\end{subfigure}
  \begin{subfigure}{.29\textwidth}
		  \includegraphics[width=1\linewidth]{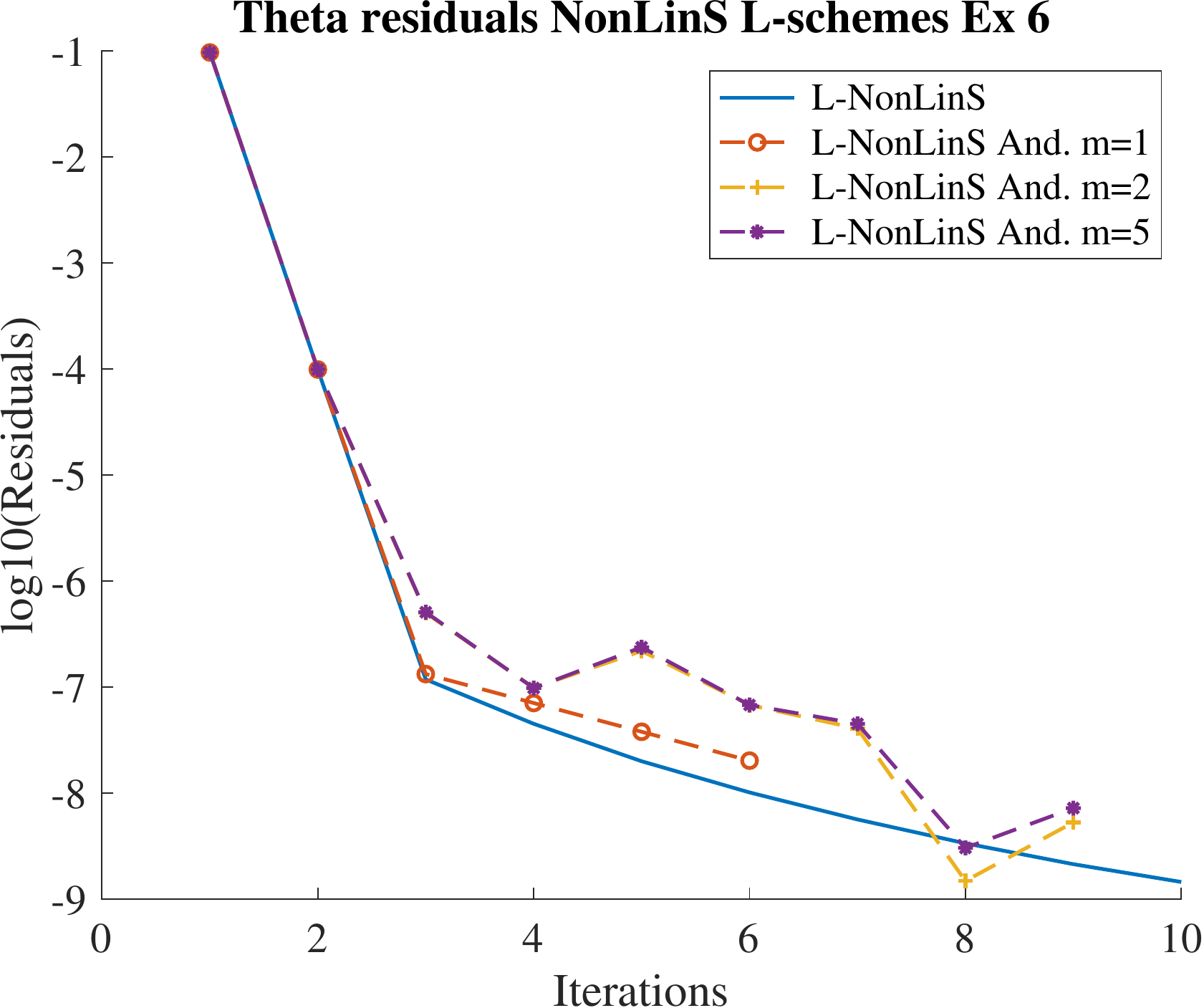}
			\caption{Water content residuals.}
			\end{subfigure}
			  \begin{subfigure}{.29\textwidth}
\includegraphics[width=1\linewidth]{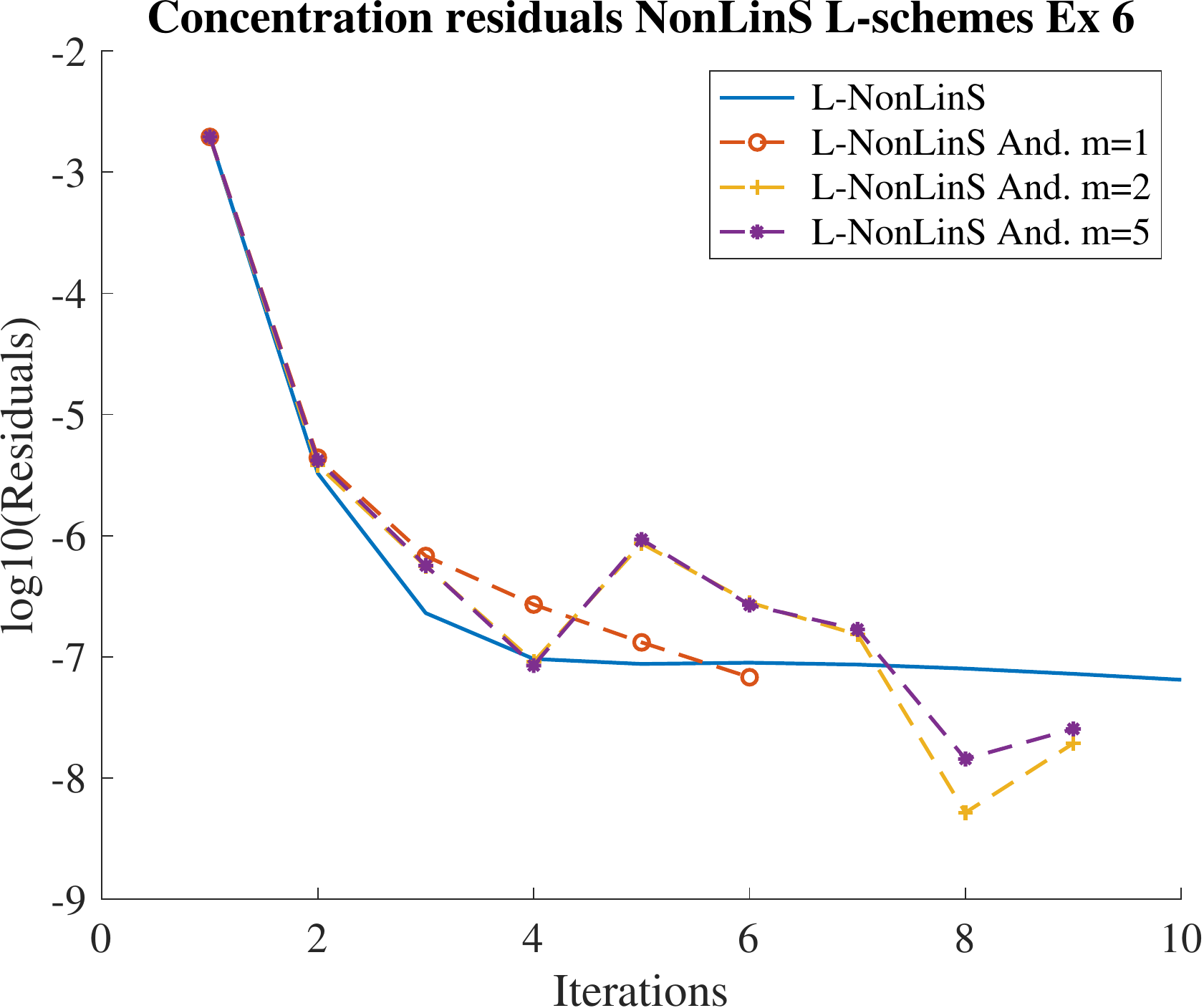}
			\caption{Concentration residuals.}
			\end{subfigure}		
		\caption{Example 5: Residuals of each unknown at the final time step, NonLinS L-scheme. Here, $L_1=L_2=L_3=0.1$, different $m$ are tested , $dx=1/40$, and $\Delta t=T/25$.}
		\label{fig:Ex5ResidualsNonLinS}
	\end{center}
\end{figure}

Finally, in the Figures \ref{fig:Ex5ResidualsMono}, \ref{fig:Ex5ResidualsAltLinS} and \ref{fig:Ex5ResidualsNonLinS} we report the residuals at the final time step, for each unknown and each algorithm, and for different values of $m$. The results are coherent with the ones presented in the Tables \ref{tab:5.0dx} and \ref{tab:5.0dt}. The monolithic solver shows a clear improvement thanks to the AA, for both $m=1$ and $m=2$. For the alternate linearized splitting solver, the AA does not seem to produce any improvement, as the rates of convergence seem to be worse. The result does not directly contradict the ones presented in the Tables \ref{tab:5.0dx} and \ref{tab:5.0dt}; it simply states that, at the final time step, the AA does not produce any improvement. On the full simulation, we could observe a clear reduction in the numbers of iterations.
For the nonlinear splitting, we can observe some improvements but once more they are not as evident as for the monolithic solvers. We already observed that both splitting solvers became unstable once the AA was applied, either requiring a larger L or a smaller $m$.

Table \ref{tab:5.1EOC} presents the precise rates of convergence of the different linearization schemes. Once more we can observe as the AA improves the rates of convergence of the solvers based on the L-scheme.

\begin{table}[h!]
\begin{center}
	\captionsetup{justification=centering,margin=2cm}
\begin{tabular}{|l|l|l|l|l|l|l|}
\hline
&\scriptsize{LS-Mono} & \scriptsize{LS-Mono And.} & \scriptsize{LS-NonLinS} & \scriptsize{LS-NonLinS And.} & \scriptsize{LS-AltLinS} & \scriptsize{LS-AltLinS And.}\\ \hline
$\Psi$   & 0.94   & 1.30 &  0.98   &  1.87  & 0.98 &  1.21    \\ \hline
$c$        & 1.01   &  1.78  &  0.95   &  1.42  & 0.93 &  1.80 \\ \hline
$\theta$ & 0.95 & 1.34   &  0.99   &  1.35  & 0.95 &  1.02    \\ \hline
\end{tabular}
\caption{Example 5: Order of convergence  of the linearization schemes.}
\label{tab:5.1EOC}
\end{center}
\end{table}

\section{Conclusions}
\label{conclusion}
We consider models for flow and reactive transport in a porous medium. Next, to account for the influence of the solute concentration on the flow parameters, we incorporate effects like dynamic capillary pressure and hysteresis. The problem results being fully coupled.

For solving the time discrete equations \eqref{richards1}, obtained after applying the Euler implicit scheme, we investigate different approaches: a monolithic solution algorithm and two splitting ones. Furthermore, for solving the nonlinear problem, two linearizations are studied: the Newton method and the L-scheme. The latter appears to be more stable than the former, which is more commonly implemented.

Finally, we have studied the effects of the Anderson acceleration. We observed that its implementation is particularly simple and can result in significant improvements. There were cases in which the differences between the accelerated and non-accelerated schemes were minimal, but due to its simplicity and the possibility of the great reduction in the numbers of iterations, we think it should always be tested. Particularly, one can either invest time in finding the optimal L parameters or the best depth $m$ for which the AA results in the fastest scheme. Often, finding the most suitable $m$ is simpler, and it can results in impressive improvements. 

\begin{acknowledgments}
	The research of D. Illiano was funded by VISTA, a collaboration between the NorwegianAcademy of Science and Letters and Equinor, project number 6367, project name: adaptive model and solver
	 simulation of enhanced oil recovery. The research of J.W. Both was supported by the Research Council of Norway
	 Project 250223, as well as the FracFlow project funded by Equinor through Akademiaavtalen. The research of I.S.
	 Pop was supported by the Research Foundation-Flanders (FWO), Belgium through the Odysseus programme (project G0G1316N) and Equinor through the Akademia grant.
\end{acknowledgments}

\end{document}